\documentclass{article}
\usepackage{theorem,amsmath,amssymb}
\newtheorem{Th}{Theorem}[section]
\newtheorem{Le}{Lemma}[section]
\newtheorem{Co}{Corollary}[section]

\newtheorem{Rem}{Remark}[section]

\date{}
\title{On the behavior of solutions of the nonstationary Stokes system near the vertex of a cone}
\author{Vladimir Kozlov and J\"urgen Rossmann}
\begin{document}
\maketitle

\abstract{The paper deals with the Dirichlet problem for the nonstationary Stokes system in a threedimensional cone.
The authors study the asymptotics of the solutions  near the vertex of the cone. \\


}

\section*{Introduction}

The present paper deals with the Dirichlet problem for the nonstationary Stokes system
\begin{eqnarray} \label{stokes1}
&&\frac{\partial u}{\partial t} - \Delta u + \nabla p = f, \quad -\nabla\cdot u = g \ \mbox{ in } K \times (0,\infty), \\ \label{stokes2}
&& u(x,t)=0 \ \mbox{ for }x\in \partial K, \ t>0, \quad u(x,0)=0 \ \mbox{ for }x\in K,
\end{eqnarray}
where $K=\{ x\in {\Bbb R}^3: \, |x|^{-1}x\in \Omega\}$ is a three-dimensional cone with vertex at the origin.
Here, $\Omega$ is a domain on the unit sphere with boundary $\partial\Omega$ of the
class $C^{2,\alpha}$, where $\alpha$ is an arbitrarily small positive number.
In our recent paper \cite{k/r-16}, we proved existence, uniqueness and regulartiy assertions for solutions
of this problem  in weighted Sobolev spaces. Analogous results for the two-dimensional case were obtained in \cite{r-18}.
The goal of the present paper is to describe the asymptotics of the solutions near the vertex of the cone $K$.
The behavior of solutions of elliptic boundary value problems (including the stationary Stokes system) near conical points or edges
is well-studied in many papers (see e.~g. the bibliography in the monographs \cite{kmr-97,kmr-01,mr-10,np-94}).
For the stationary Stokes system, we refer also to \cite{Dauge-89,Kozlov/M/S-94,mp-83,mps-79,sol-79}.
Furthermore, there are several papers dealing with the behavior of solutions of the heat equation near conical points and edges
(see \cite{dCN-11,k/m-87,k/r-11,kweon,nazarov-01,Sol-84,Sol-01}). In particular, it was shown in \cite{k/m-87} that the solution
of the first initial boundary value problem for the heat equation in a 3-dimensional cone is a sum of ``singular" terms
\[
r^{\lambda_j}\, \phi_j(\omega)\, (r^2\partial_t)^k \, ({\cal E}h_j)(r,t)
\]
and a regular remainder. Here, $r=|x|$, $\omega=|x|^{-1}x$,  $\phi_j$ is an eigenfunction of the Beltrami operator
$-\delta$ on $\Omega$ corresponding to the eigenvalue $\Lambda_j=\lambda_j(1+\lambda_j)$, $h_j$ is a function in $t$ and ${\cal E}$
is an extension operator. A similar result was obtained in \cite{kozlov-88,kozlov-89,kozlov-91} for a class of general parabolic problems.
However, the class of problems considered in \cite{kozlov-88,kozlov-89,kozlov-91} does not include the Stokes system.
In fact, in a number of results, there are essential differences between the Stokes system on one hand and the heat equation
(and other parabolic problems) on the other hand. For all these problem, one has existence and uniqueness results in weighted Sobolev
spaces. For the heat equation, and also for the class of problems in \cite{kozlov-88,kozlov-89,kozlov-91}, the bounds for
the weight parameter depend only on the eigenvalues of one operator pencil. In the case of the Stokes system, the eigenvalues
of two different operator pencils appear in these bounds. Concerning the behavior of the solutions at infinity, we have
completely different assertions for the heat equation and the Stokes system (see \cite{k/r-18a}).
Nevertheless, similarly as in \cite{kozlov-88,kozlov-89,kozlov-91}, we can split the solution of the problem
(\ref{stokes1}), (\ref{stokes2}) into finitely many singular terms and a regular remainder. The singular terms depend on the eigenvalues $\lambda_j$
of the operator pencil generated by the stationary Stokes system in a cone.

The major part of the paper (Section 1 and 2) deals with the parameter-depending problem
\begin{equation} \label{par1}
s\, \tilde{u} - \Delta \tilde{u} + \nabla \tilde{p} = \tilde{f}, \quad -\nabla\cdot \tilde{u} = \tilde{g}  \
   \mbox{ in } K, \quad \tilde{u} = 0 \ \mbox{ on } \partial K,
\end{equation}
which arises after the Laplace transformation with respect to the time $t$, for $\mbox{Re}\, s \ge 0$, $s\not=0$.
In Section 1, we present an existence and uniqueness theorem for this problem which was proved in \cite{k/r-16}.
Let $V_\beta^l(K)$ be the weighted Sobolev space of all functions (vector-functions) with finite norm
\begin{equation} \label{Vbeta}
\| u\|_{V_\beta^l(K)} = \Big( \int_K \sum_{|\alpha|\le l} r^{2(\beta-l+|\alpha|)}\big| \partial_x^\alpha u(x)\big|^2\, dx\Big)^{1/2},
\end{equation}
while $E_\beta^l(K)$ is the  weighted Sobolev space with the norm
\begin{equation} \label{Ebeta}
\| u\|_{E_\beta^l(K)} = \Big( \int_K \sum_{|\alpha|\le l} \big( r^{2\beta}+r^{2(\beta-l+|\alpha|)}\big)\big|
  \partial_x^\alpha u(x)\big|^2\, dx\Big)^{1/2},
\end{equation}
$r=r(x)$ denotes the distance of the point $x$ from the vertex of the cone.
By Theorem \ref{t1}, there exists a uniquely determined solution $(\tilde{u},\tilde{p}) \in E_\beta^2(K)\times V_\beta^1(K)$
for arbitrary $\tilde{f}\in E_\beta^0(K)$ and $\tilde{g}\in V_\beta^1(K)\cap(V_{-\beta}^1(K))^*$ if the weight parameter $\beta$
satisfies the inequalities
\begin{equation} \label{interval}
\frac 12 - \lambda_1 <\beta  < \min\Big(\mu_2 +\frac 12\, , \lambda_1 +\frac 32\Big), \ \ \beta\not=\frac 12
\end{equation}
and $\int_K \tilde{g}(x)\, dx=0$ for $\beta>\frac 12$. Here, $\lambda_1$ and $\mu_2$ are positive numbers depending on the cone.
More precisely, $\lambda_1$ is the smallest positive eigenvalue of the operator pencil ${\cal L}(\lambda)$ generated by the Dirichlet problem for
the stationary Stokes system, while $\mu_2$ is the smallest positive eigenvalue of the operator pencil
${\cal N}(\lambda)$ generated by the Neumann problem for the Laplacian, respectively (see Section 1).
Since $\lambda=1$ is always an eigenvalue of the pencil ${\cal L}(\lambda)$, the eigenvalue $\lambda_1$ is not greater than 1.

Section 2 is concerned with the asymptotics of the solutions $(\tilde{u},\tilde{p})$ of the problem (\ref{par1}) near the vertex of the cone. Suppose that
$\tilde{u} \in E_\beta^2(K)$, $\tilde{p}\in V_\beta^1(K)$, where $\beta$ satisfies the inequalities (\ref{interval}), and that the data
$\tilde{f}$, $\tilde{g}$ belong to the spaces $E_\gamma^0(K)$ and $V_\gamma^1(K)\cap (V_{-\gamma}^1(K))^*$, respectively, where $\gamma < \beta$.
We show in Section 2, that $(\tilde{u},\tilde{p})$ is a sum of terms
\[
\zeta\big(|s|r^2\big)\, c_{j,k}(s)\, \big( u_{j,k,\gamma}(x,s),p_{j,k,\gamma}(x,s)\big)
\]
and a remainder $(\tilde{v},\tilde{q}) \in E_\gamma^2(K)\times V_\gamma^1(K)$.
Here, $\zeta$ is a cut-off function on $(0,\infty)$ with support in $[0,1]$, $\zeta(r)=1$ for $r<1/2$, $c_{j,k}$ are constants depending on $s$ and on the data,
and $u_{j,k,\gamma}$, $p_{j,k,\gamma}$ are functions of the form
\[
u_{j,k,\gamma}(x,s)= r^{\lambda_j}\sum_\mu (sr^2)^\mu {\cal P}_\mu(\omega,\log r), \quad
p_{j,k,\gamma}(x,s)= r^{\lambda_j-1}\sum_\mu (sr^2)^\mu {\cal Q}_\mu(\omega,\log r),
\]
where $\lambda_j$ are eigenvalues of the pencil ${\cal L}(\lambda)$, ${\cal P}_\mu(\omega,\log r)$, ${\cal Q}_\mu(\omega,\log r)$
are polynomials in $\log r$ with coefficients depending on $\omega=x/|x|$. Furthermore, we prove a formula for the coefficients $c_{j,k}$
and obtain an upper bound for $|c_{j,k}|$ depending on $s$ (Theorem \ref{At1}). In the case $\gamma > -\mu_2-\frac 12$, we get the estimate
\[
\| \tilde{v}\|^2_{V_\gamma^2(K)} + |s|^2\, \| \tilde{v}\|^2_{V_\gamma^0(K)}
   \le c\, \Big( \| \tilde{f}\|^2_{V_\gamma^0(K)} +  \| \tilde{g}\|^2_{V_\gamma^1(K)}+ |s|^2\, \| \tilde{g}\|^2_{(V_{-\gamma}^1(K))^*}\Big)
\]
for the remainder $\tilde{v}$ (Theorem \ref{At2}). The analogous estimate for the $V_\gamma^1(K)$-norm of $\tilde{q}$ holds only if $\gamma>-\frac 12$.
In the last case, we have the simpler representation
\[
u_{j,k,\gamma}(x,s)= u_{j,k}^{(0)}(x) = r^{\lambda_j}\, \phi_{j,k}(\omega), \quad p_{j,k,\gamma}(x,s)= p_{j,k}^{(0)}(x) = r^{\lambda_j-1}\, \psi_{j,k}(\omega)
\]
for the singular terms.

The last section of the paper deals with solutions of the time-dependent problem in the weighted space $W_\beta^{2,1}(K\times{\Bbb R}_+)\times
L_2({\Bbb R}_+,V_\beta^1(K))$. Here, $W_\beta^{2,1}(K\times{\Bbb R}_+)$ is the space of all $u\in L_2({\Bbb R}_+,V_\beta^2(K))$ such that
$\partial_t u \in L_2({\Bbb R}_+,V_\beta^0(K))$. By \cite[Theorem 3.1]{k/r-16}, the problem (\ref{stokes1}), (\ref{stokes2}) has a uniquely
determined solution in this space if $f\in L_2({\Bbb R}_+,V_\beta^0(K))$, $g\in L_2({\Bbb R}_+,V_\beta^2(K))$,
$\partial_t g\in L_2({\Bbb R}_+,V_\beta^1(K))$, $\beta$ satisfies the inequalities (\ref{interval}) and $g(x,0)=0$ (in addition $g$ must satisfy the condition $\int_K g(x,t)\, dx=0$ if $\beta>\frac 12$).
Applying the results of Section 2, we obtain the asymptotics of this solution at the vertex of the cone if
$f \in L_2({\Bbb R}_+,V_\gamma^0(K))$, $g\in L_2({\Bbb R}_+,V_\gamma^1(K))$, $\partial_t g \in L_2({\Bbb R}_+,(V_{-\gamma}^1(K))^*)$ and $\gamma <\beta$.
For example, we obtain the decomposition
\[
(u,p) = \sum_{0<\lambda_j<-\gamma+1/2} \sum_{k=1}^{\kappa_j}  H_{j,k}(r,t) \, \big( u_{j,k}^{(0)}(x),p_{j,k}^{(0)}(x)\big)+ (v,q)
\]
with a remainder $(v,q)\in W_\gamma^{2,1}(K\times{\Bbb R}_+)\times L_2({\Bbb R}_+,V_\gamma^1(K))$ if $-\frac 12 < \gamma < \frac 12 - \lambda_1$ (see Theorem \ref{t8}).
For the velocity $u$, the decomposition 
\[
u = \sum_{j,k}  u_{j,k,\gamma}(x,\partial_t)\, H_{j,k}(r,t) + v
\]
with a remainder $v\in W_\gamma^{2,1}(Q)$ holds even if $\gamma>-\mu_2-\frac 12$. The coefficients $H_{j,k}$ are given by the formula
\begin{equation} \label{7t8}
H_{j,k}(r,t) = -\int_0^t \int_K \big( f(y,\tau)\, K_{j,k}(x,y,t-\tau) + g(y,\tau)\, T_{j,k}(x,y,t-\tau)\big)\, dy\, d\tau
\end{equation}
and are extensions of certain functions $h_{j,k}(t)$. The kernels $K_{j,k}$ and $H_{j,k}$ in (\ref{7t8}) satisfy the estimates of Lemma \ref{l2}.
Under additional assumptions on the $t$-derivatives of $f$ and $g$, we obtain a decomposition of the solution with a remainder 
$(v,q)\in W_\gamma^{2,1}(K\times{\Bbb R}_+)\times L_2({\Bbb R}_+,V_\gamma^1(K))$, $\gamma<-\frac 12$ (see Theorem \ref{t8a}).

\section{Solvability of the parameter-depending problem}

Let $s$ be an arbitrary complex number, $\mbox{Re}\, s \ge 0$. We consider the boundary value problem
\begin{equation} \label{par3}
s\, u - \Delta u + \nabla p = f, \ \ - \nabla\cdot u = g \ \mbox{ in }K, \quad u=0\ \mbox{ on }\partial K\backslash \{ 0\}.
\end{equation}
For nonnegative integer $l$ and real $\beta$, we define the weighted Sobolev spaces $V_\beta^l(K)$ and $E_\beta^l(K)$
as the sets of all functions (or vector functions) with finite norms (\ref{Vbeta}) and (\ref{Ebeta}), respectively.
Note that the spaces $V_\beta^l(K)$ and $E_\beta^l(K)$ can be also defined as the closures of
$C_0^\infty(\overline{K}\backslash \{ 0\})$ with respect to the above norms.
Furthermore, we define $\stackrel{\circ}{V}\!{}_\beta^1(K)$ and $\stackrel{\circ}{E}\!{}_\beta^1(K)$
as the spaces of all functions $u\in V_\beta^1(K)$ and $u\in E_\beta^1(K)$, respectively, which are zero on
$\partial K\backslash \{ 0\}$. Moreover, let
\[
X_\beta^1(K) = V_\beta^1(K) \cap \big( V_{-\beta}^1(K)\big)^*.
\]
It can be easily shown that the operator of the problem (\ref{par3}) realizes a continuous mapping
\begin{eqnarray*}
&& \big( E_\beta^2(K) \cap \stackrel{\circ}{E}\!{}_{\beta}^1(K)\big) \times V_\beta^1(K) \ni (u,p) \\
&& \qquad \to \  \big( su-\Delta u + \nabla p, - \nabla\cdot u\big) \in E_\beta^0(K)\times X_\beta^1(K).
\end{eqnarray*}
for arbitrary $s$. We denote this operator by $A_\beta$.

We introduce the following operator pencils ${\cal L}(\lambda)$ and ${\cal N}(\lambda)$ generated by
the Dirichlet problem for the stationary Stokes system and the Neumann problem for the Laplacian
in the cone $K$, respectively. For every complex $\lambda$, we define the operator ${\cal L}(\lambda)$ as the mapping
\begin{eqnarray*}
&& \stackrel{\circ}{W}\!{}^1(\Omega)\times L_2(\Omega) \ni \left( \begin{array}{c} U \\ P \end{array}\right) \\ && \qquad
  \to \left( \begin{array}{c} r^{2-\lambda}\big(-\Delta r^{\lambda}U(\omega)+\nabla r^{\lambda-1}P(\omega)\big)\\*[1ex]
  -r^{1-\lambda}\nabla\cdot\big( r^\lambda U(\omega)\big)\end{array}\right) \in W^{-1}(\Omega)\times L_2(\Omega),
\end{eqnarray*}
where $r=|x|$ and $\omega=x/|x|$. The properties of the pencil ${\cal L}$ are studied, e.g., in \cite{kmr-01}.
In particular, it is known that the numbers $\lambda$, $\bar{\lambda}$ and $-1-\lambda$ are simultaneously eigenvalues
of the pencil ${\cal L}(\lambda)$ or not. The eigenvalues in the strip $-2\le \mbox{Re}\, \lambda \le 1$ are real,
and the numbers $1$ and $-2$ are always eigenvalues of the pencil ${\cal L}(\lambda)$. If $\Omega$ is contained in a half-sphere,
then $\lambda=1$ and $\lambda=-2$ are the only eigenvalues in the interval $[-2,1]$ (cf. \cite[Theorem 5.5.5]{kmr-01}).
We denote by $\lambda_j$ the eigenvalues with positive real part and by $\lambda_{-j}=-1-\lambda_j$
the eigenvalues with negative real part, $j=1,2,\ldots$, and state that
\[
\cdots \le \mbox{Re}\, \lambda_{-2} < \lambda_{-1} <  \lambda_1 < \mbox{Re}\, \lambda_2 \le \cdots
\]
Here, $0<\lambda_1\le 1$ and $-2\le \lambda_{-1}<-1$.
Note that the eigenvalues $\lambda_j$ and $\lambda_{-j}$ have the same geometric and algebraic multiplicities.
The operator ${\cal N}(\lambda)$ is defined as
\[
{\cal N}(\lambda)\, U = \Big( -\delta U-\lambda(\lambda+1)U\, , \, \frac{\partial U}{\partial\vec{n}}\big|_{\partial\Omega}\Big)\
  \quad\mbox{for } U\in W^2(\Omega).
\]
As is known (see e.g. \cite[Section 2.3]{kmr-01}), the eigenvalues of this pencil are real and
generalized eigenfunctions do not exist. The spectrum contains, in particular, the simple eigenvalues $\mu_1=0$ and $\mu_{-1}=-1$ with the
eigenfunction $\phi_1 = const.$ The interval $(-1,0)$ is free of eigenvalues. Let $\mu_j$ be the nonnegative
and $\mu_{-j}=-1-\mu_j$ the negative eigenvalues of the pencil ${\cal N}(\lambda)$, $j=1,2,\ldots$, where
\[
\cdots < \mu_{-2} < -1 =\mu_{-1} < \mu_1 =0 < \mu_2 < \cdots.
\]
Obviously, $\mu_j$ and $\mu_{-j}$ are the solutions of the equation $\mu(\mu+1)=-M_j$, where $M_j$ is the $j$th eigenvalue
of the operator $-\delta$ with Neumann boundary condition. As was shown in \cite[Theorem 2.1]{k/r-16}, the operator
$A_\beta$ has closed range and finite-dimensional kernel if the line $\mbox{Re}\, \lambda = -\beta+1/2$ does not contain
eigenvalues of the pencil ${\cal L}(\lambda)$ and $-\beta-1/2$ is not an eigenvalue of the pencil ${\cal N}(\lambda)$.
Furthermore, the following result was proved in \cite[Lemma 2.10, Theorem 2.4]{k/r-16}

\begin{Th} \label{t1}
Suppose that $\mbox{\em Re}\, s \ge 0$ and $s\not=0$.

{\em 1)} If $-\mu_2 -1/2 < \beta < \lambda_1 +3/2$, then $A_\beta$ is injective.

{\em 2)} If $\frac 12 -\lambda_1 <\beta < \frac 12$, then the operator $A_\beta$ is an isomorphism onto $E_\beta^0(K) \times X_\beta^1(K)$,
and the estimate
\begin{equation} \label{1t1}
\| u\|_{V_\beta^2(K)} + |s|\, \| u\|_{V_\beta^0(K)} + \| p\|_{V_\beta^1(K)} \le c\, \Big( \| f\|_{V_\beta^0(K)}
  + \| g \|_{V_\beta^1(K)} + |s|\, \| g\|_{(V_{-\beta}^1(K))^*}\Big)
\end{equation}
is valid for every solution $(u,p)\in E_\beta^2(K)\times V_\beta^1(K)$ of the problem {\em (\ref{par3})}. Here, the constant $c$ is independent
of $f,g$ and $s$.

{\em 3)} If $\frac 12 <\beta < \min\big( \mu_2 +\frac 12\, , \, \lambda_1 +\frac 32\big)$, then
the operator $(u,p) \to (f,g)$ of the problem {\em (\ref{par3})} is an isomorphism from
$\big( E_\beta^2(K) \cap\stackrel{\circ}{E}\!{}_\beta^1(K)\big)\times V_\beta^1(K)$ onto the subspace
of all $(f,g) \in E_\beta^0(K) \times X_\beta^1(K)$ satisfying the condition
\begin{equation} \label{2t1}
\int_K g(x)\, dx =0.
\end{equation}
Furthermore, the estimate {\em (\ref{1t1})} is valid for every solution $(u,p)\in E_\beta^2(K)\times V_\beta^1(K)$ of the problem {\em (\ref{par3})}.
\end{Th}

The conditions on $\beta$ in Theorem \ref{t1} are sharp if $\lambda_1<1$.  In the case that $\lambda_1=1$ and $\lambda_1$ is a simple eigenvalue
of the pencil ${\cal L}(\lambda)$, an analogous  existence and uniqueness result holds even for
$\max\big(-\mu_2-\frac 12,\, \frac 12-\mbox{Re}\, \lambda_2\big) < \beta < \min\big(\mu_2+\frac 12,\, \mbox{Re}\, \lambda_2+\frac 32\big)$,
$\beta\not= \pm \frac 12$, $\beta\not= \frac 52$. Here $g$ must satisfy the condition (\ref{2t1}) if $\frac 12 < \beta <\frac 52$
(see \cite{k/r-18a}).

\begin{Rem} \label{r1a}
{\em For $s=0$ (i.e., for the stationary Stokes system), the mapping
\[
(u,p) \to (-\Delta u+\nabla p,-\nabla\cdot u)
\]
is continuous from $(V_\beta^2(K)\cap \stackrel{\circ}{V}\!{}_\beta^1(K))\times V_\beta^1(K)$ into $V_\beta^0(K)\times V_\beta^1(K)$.
This mapping is bijective if and only if the line $\mbox{Re}\, \lambda=-\beta+1/2$ is free of eigenvalues of the pencil ${\cal L}(\lambda)$.
Under this condition every solution $(u,p)\in V_\beta^2(K) \times V_\beta^1(K)$ of the problem (\ref{par3}) also satisfies the inequality
(\ref{1t1}) (see, e.~g., \cite{mp-78a}, \cite[Chapter 3, Theorem 5.1]{np-94}).}
\end{Rem}

Furthermore, the following regularity result will be used in this paper.

\begin{Le}  \label{l5b}
Suppose that $(u,p) \in E_\beta^2(K) \times V_\beta^1(K)$ is a solution of the problem {\em (\ref{par3})}, where $\mbox{\em Re}\, s\ge 0$
and $|s|=1$, $f\in E_\beta^0(K) \cap E_\gamma^0(K)$ and $g\in X_\beta^1(K)\cap X_\gamma^1(K)$.
We assume that one of the following two conditions is satisfied:
\begin{itemize}
\item[{\em (i)}] $\beta<\gamma$ and the interval $-\gamma-1/2 \le \lambda \le -\beta-1/2$ does not contain eigenvalues of the pencil
${\cal N}(\lambda)$,
\item[{\em (ii)}] $\beta>\gamma$ and the strip $-\beta+1/2 \le \mbox{\em Re}\, \lambda\le -\gamma+1/2$ is free of eigenvalues
  of the pencil ${\cal L}(\lambda)$.
\end{itemize}
Then $u\in E_\gamma^2(K)$, $p\in V_\gamma^1(K)$ and
\[
\| u\|_{E_\gamma^2(K)} +\| p\|_{V_\gamma^1(K)} \le c\, \Big( \| f\|_{E_\gamma^0(K)} + \| g\|_{X_\gamma^1(K)}
  + \| u\|_{E_\beta^2(K)} +\| p\|_{V_\beta^1(K)}\Big).
\]
Here, the constant $c$ is independent of $f$, $g$ and $s$.
\end{Le}

\section{Asymptotics of solutions of the parameter-depending problem}

We consider a solution $(u,p) \in E_\beta^2(K)\times V_\beta^1(K)$ of the problem (\ref{par3}) with data
$f\in E_\beta^0(K)\cap E_\gamma^0(K)$ and $g\in X_\beta^1(K)\cap X_\gamma^1(K)$, where $\gamma<\beta$. If the condition of Lemma \ref{l5b}
is not satisfied, i.~e., if the strip $-\beta+1/2 \le \mbox{Re}\, \lambda\le -\gamma+1/2$ contains eigenvalues
of the pencil ${\cal L}(\lambda)$, we cannot expect that $u\in E_\gamma^2(K)$ and $p\in E_\gamma^1(K)$. We show that then the solution
is a finite sum of ``singular'' terms and a remainder $(v,q)\in E_\gamma^2(K)\times V_\gamma^1(K)$.

\subsection{Singularities of solutions near the vertex of the cone}

As was mentioned above, the eigenvalues $\lambda_i$ and $\lambda_{-i}$ have the same algebraic multiplicity $\kappa_i=\kappa_{-i}$.
Thus for every eigenvalue $\lambda_i$ (with positive or negative $i$), there exist $\kappa_i$ linearly independent solutions
$(u_{i,k}^{(0)},p_{i,k}^{(0)})$, $k=1,2,\ldots,\kappa_i$ of the problem
\begin{equation} \label{hom}
-\Delta u + \nabla p =0, \ \ \nabla\cdot u=0\ \mbox{ in }K, \quad u=0\ \mbox{ on }\partial K\backslash \{ 0\}
\end{equation}
which have form
\begin{equation} \label{sing}
u_{i,k}^{(0)}(x) = r^{\lambda_i} \sum_{\nu=0}^\sigma \frac{(\log r)^\nu}{\nu !}\, \phi_{i,k,\sigma-\nu}(\omega), \quad
p_{i,k}^{(0)}(x) = r^{\lambda_i-1} \sum_{\nu=0}^\sigma \frac{(\log r)^\nu}{\nu !}\, \psi_{i,k,\sigma-\nu}(\omega), \quad
\end{equation}
where $(\phi_{i,k,\nu},\psi_{j,k,\nu})$ are eigenvectors (if $\nu=0$) or generalized eigenvectors (if $\nu>0$) of the
pencil ${\cal L}(\lambda)$ with respect to the eigenvalue $\lambda_j$. The pairs $(u_{i,k}^{(0)},p_{i,k}^{(0)})$, $k=1,\ldots,\kappa_i$
form a basis in the space of the solutions of the problem (\ref{hom}) which have the form (\ref{sing}).

Now let $i\ge 1$, and let $\gamma$ be an arbitrary real number, $\gamma<\frac 12-\lambda_i$. We define $N_{i,\gamma}$ as the smallest integer $N$
such that $\gamma +\mbox{Re}\, \lambda_j + 2N > - \frac 32$. Then for $\mu=1,2,\ldots,N_{i,\gamma}$ there exist vector functions
$(u_{i,k}^{(\mu)},p_{i,k}^{(\mu)})$ of the form
\begin{equation} \label{uknu1}
u_{i,k}^{(\mu)}(x) = r^{\lambda_i} \sum_{\mu=0}^\sigma (\log r)^\mu\ a_\mu(\omega), \quad
p_{j,k}^{(0)}(x) = r^{\lambda_i-1} \sum_{\mu=0}^\sigma (\log r)^\mu \, b_\mu(\omega),
\end{equation}
such that
\begin{eqnarray} \label{uknu2}
&& - \Delta (r^{2\mu} u_{i,k}^{(\mu)}) + \nabla (r^{2\mu}p_{i,k}^{(\mu)}) = -r^{2\mu-2}u_{i,k}^{(\mu-1)}, \quad \nabla\cdot (r^{2\mu}u_{i,k}^{(\mu)})=0
  \ \mbox{ in }K, \\ \label{uknu3}
&& u_{i,k}^{(\mu)}=0\ \mbox{ on }\partial K\backslash \{ 0\}
\end{eqnarray}
(cf. \cite[Lemma 5.1.6]{kmr-97}). Note that the pair $(u_{i,k}^{(\mu)},p_{i,k}^{(\mu)})$ in (\ref{uknu2}), (\ref{uknu3}) is uniquely determined
for a given vector function $u_{i,k}^{(\mu-1)}$ if $\lambda_i+2\mu$ is not an eigenvalue of the pencil ${\cal L}(\lambda)$. If however
$\lambda_i+2\mu=\lambda_m$ is an eigenvalue, then, together with the vector function $(u_{i,k}^{(\mu)},p_{i,k}^{(\mu)})$, the vector functions
\[
(u_{i,k}^{(\mu)},p_{i,k}^{(\mu)}) + r^{-2\mu} \sum_{\sigma=1}^{\kappa_m} c_\sigma \, (u_{m,\sigma}^{(0)},p_{m,\sigma}^{(0)})
\]
are also solutions of the problem (\ref{uknu2}), (\ref{uknu3}) which have the form (\ref{uknu1}). We define
\begin{equation} \label{Uikg}
u_{i,k,\gamma}(x,s) = \sum_{\mu=0}^{N_{i,\gamma}} (sr^2)^\mu \, u_{i,k}^{(\mu)}(x), \quad
p_{i,k,\gamma}(x,s) = \sum_{\mu=0}^{N_{i,\gamma}} (sr^2)^\mu \, p_{i,k}^{(\mu)}(x).
\end{equation}
Obviously, $u_{i,k,\gamma}=0$ on $\partial K \backslash \{ 0\}$.

\begin{Rem} \label{Ar1}
{\em The number $\lambda=1$ is always an eigenvalue of the pencil ${\cal L}(\lambda)$ with the constant eigenfunction
$(\phi,\psi)=(0,c)$. In the case $\lambda_i=1$, $u_{i,k}^{(0)}=0$, $p_{i,k}^{(0)}=c$, we can obviously set $u_{i,k}^{(\mu)}=0$, $p_{i,k}^{(\mu)}=0$
for $\mu\ge 1$, i.e. $u_{i,k,\gamma}=0$, $p_{i,k,\gamma}=p_{i,k}^{(0)}=c$ for arbitrary $\gamma< -\frac 12$.}
\end{Rem}

In the following, let $\zeta$ be a two times continuously differentiable function on $(0,\infty)$ such that $\zeta(r)=1$ for $r< 1/2$ and $\zeta(r)=0$ for $r>1$.
Furthermore, let
\[
\zeta_s(r) = \zeta\big( |s|\, r^2\big).
\]

\begin{Le} \label{Al1}
Let $i\ge 1$ and $\gamma<\frac 12 -\lambda_1$. Then
\begin{equation} \label{1Al1}
(s-\Delta)\, u_{i,k,\gamma} + \nabla p_{i,k,\gamma} = s\, (sr^2)^{N_{i,\gamma}}\, u_{i,k}^{(N_{i,\gamma})}, \quad \nabla\cdot u_{i,k,\gamma}=0
\end{equation}
in $K$. In particular, $\zeta_s\, \big( (s-\Delta)\, u_{i,k,\gamma} + \nabla p_{i,k,\gamma}\big) \in V_\gamma^0(K)$
\end{Le}

P r o o f. From the definition of $u_{i,k,\gamma}$ and $p_{i,k,\gamma}$ it follows that
\begin{eqnarray*}
&& (s-\Delta)\, u_{i,k,\gamma} + \nabla p_{i,k,\gamma} = s\sum_{\mu=0}^{N_{i,\gamma}} (sr^2)^\mu \, u_{i,k}^{(\mu)} +
  \sum_{\mu=0}^{N_{i,\gamma}} s^{\mu}\Big( -\Delta (r^{2\mu} \, u_{i,k}^{(\mu)})+\nabla(r^{2\mu}p_{i,k}^{(\mu)})\Big)  \\
&& = s\sum_{\mu=0}^{N_{i,\gamma}} (sr^2)^\mu \, u_{i,k}^{(\mu)} - s\sum_{\mu=1}^{N_{i,\gamma}} (sr^2)^{\mu-1} \, u_{i,k}^{(\mu-1)}
  = s\, (sr^2)^{N_{i,\gamma}}\, u_{i,k}^{(N_{i,\gamma})}.
\end{eqnarray*}
Obviously, $\zeta_s\, r^{2N_{i,\gamma}}\, u_{i,k}^{(N_{i,\gamma})}\in V_\gamma^0(K)$ if $\gamma+2N_{i,\gamma}+\mbox{Re}\, \lambda_i > -3/2$. This
proves the lemma. \hfill $\Box$ \\

Let $(u,p)\in E_\beta^2(K)\times V_\beta^1(K)$ be a solution of the problem (\ref{par3}). We assume that $\beta<\lambda_1+\frac 32$. Then the
half-plane $\mbox{Re}\, \lambda > -\beta+1/2$ contains only eigenvalues of the pencil ${\cal L}(\lambda)$ with positive real part.

\begin{Le} \label{Al2}
Let $(u,p) \in E_\beta^2(K) \times V_\beta^1(K)$ be the solution of the problem {\em (\ref{par3})}, where
\[
f\in E_\beta^0(K)\cap E_\gamma^0(K), \quad g\in X_\beta^1(K)\cap X_\gamma^1(K), \quad \gamma < \beta < \lambda_1 +\frac 32.
\]
Suppose that the lines $\mbox{\em Re}\, \lambda = \frac 12 -\beta$ and $\mbox{\em Re}\, \lambda = \frac 12 -\gamma$
do not contain eigenvalues of the pencil ${\cal L}(\lambda)$. Then
\begin{equation} \label{1Al2}
( u,p) = \zeta_s \sum_{i\in I_{\beta,\gamma}} \sum_{k=1}^{\kappa_i} c_{i,k}(s)\, \big( u_{i,k,\gamma},p_{i,k,\gamma}\big) + (v,q),
\end{equation}
where $I_{\beta,\gamma}$ is the set of all $i$ such that $\frac 12 -\beta < \mbox{\em Re}\, \lambda_i < \frac 12 - \gamma$,
$c_{i,k}$ are constants depending on $s$, and $(v,q)\in E_\gamma^2(K)\times V_\gamma^1(K)$. If $|s|=1$, then the inequality
\begin{eqnarray} \label{3Al2} \nonumber
&&\| v\|_{E_\gamma^2(k)} + \| q\|_{V_\gamma^1(K)} +  \sum_{i\in I_{\beta,\gamma}} \sum_{k=1}^{\kappa_i} |c_{i,k}|  \\
&& \le c\, \Big( \| f\|_{E_\gamma^0(K)} + \| g\|_{V_\gamma^1(K)} + \| u\|_{E_\beta^2(K)} + \| p\|_{V_\beta^1(K)}\Big)
\end{eqnarray}
is satisfied, where $c$ is independent of $s$.
\end{Le}

P r o o f.
Note that the set $I_{\beta,\gamma}$ contains only positive integers since $\frac 12 - \beta > -1-\lambda_1=\lambda_{-1}$.

1) Suppose first that $\beta-2\le \gamma< \beta$. Then $E_\beta^2(K) \subset V_\gamma^0(K)$ and
$\gamma+\mbox{Re}\, \lambda_i\ge \beta+\mbox{Re}\, \lambda_i -2 > -\frac 32$, i.~e.,
$N_{i,\gamma}=0$ for $i\in I_{\beta,\gamma}$. This means that $u_{i,k,\gamma}=u_{i,k}^{(0)}$ and $p_{i,k,\gamma}=p_{i,k}^{(0)}$.
Since
\[
-\Delta u + \nabla p= f - s u \in V_\gamma^0(K)\quad\mbox{and}\quad -\nabla\cdot u = g \in V_\gamma^1(K),
\]
we obtain
\[
(u,p) = \sum_{i\in I_{\beta,\gamma}} \sum_{k=1}^{\kappa_i} c_{i,k}\, (u_{i,k}^{(0)},p_{i,k}^{(0)}) + (V,Q),
\]
where $V\in V_\gamma^2(K)$ and $Q\in V_\gamma^1(K)$ and
\begin{eqnarray*}
&& \| V\|_{V_\gamma^2(k)} + \| Q\|_{V_\gamma^1(K)} +  \sum_{i\in I_{\beta,\gamma}} \sum_{k=1}^{\kappa_i} |c_{i,k}| \\
&& \le c\, \Big( \| f-su\|_{V_\beta^0(K)\cap V_\gamma^0(K)} + \| g\|_{V_\beta^1(K)\cap V_\gamma^1(K)} \Big) \\
&&  \le c\, \Big( \| f\|_{V_\beta^0(K)\cap V_\gamma^0(K)} + \| g\|_{V_\beta^1(K)\cap V_\gamma^1(K)}
  +|s|\, \| u\|_{E_\beta^2(K)}\Big)
\end{eqnarray*}
(see, e.~g., \cite[Chapter 3, Theorem 5.6]{np-94}). This implies
\[
(u,p) = \zeta_s \sum_{i\in I_{\beta,\gamma}} \sum_{k=1}^{\kappa_i} c_{i,k}\, (u_{i,k}^{(0)},p_{i,k}^{(0)}) + (v,q),
\]
where
\[
(v,q) = (V,Q) + (1-\zeta_s)\, \sum_{i\in I_{\beta,\gamma}} \sum_{k=1}^{\kappa_i} c_{i,k}\, (u_{i,k}^{(0)},p_{i,k}^{(0)}) \in V_\gamma^2(K)\times V_\gamma^1(K).
\]
Obviously,
\[
(v,q) = \zeta_s\, (V,Q) + (1-\zeta_s)\, (u,p)
\]
and, consequently, $(v,q) \in E_\gamma^2(K)\times V_\gamma^1(K)$,
\begin{eqnarray*}
\| v\|_{E_\gamma^2(k)} + \| q\|_{V_\gamma^1(K)} \le
  c\, \Big( \| V \|_{V_\gamma^2(K)} + \| Q\|_{V_\gamma^1(K)} + \| u\|_{E_\beta^2(K)} + \| p\|_{V_\beta^1(K)}\Big),
\end{eqnarray*}
where $c$ depends only on $|s|$. This proves the lemma for $\beta-2\le \gamma<\beta$.

2) Suppose  that  $\beta-m\le \gamma <\beta-m+1$, where $m$ is an integer, $m\ge 3$, and that the lemma is proved for $\gamma \ge \beta-m+1$.
We can choose a number $\delta$ such that $\gamma<\beta-2\le \delta <\beta-1$ and there are no eigenvalues of the pencil
${\cal L}(\lambda)$ on the line $\mbox{Re}\, \lambda = \frac 12 -\delta$. Since $E_\delta^0(K)\subset E_\beta^0(K)\cap E_\gamma^0(K)$
and $X_\delta^1(K)\subset X_\beta^1(K)\cap X_\gamma^1(K)$, it follows from the first step of the proof that
\[
(u,p) = \zeta_s \sum_{i\in I_{\beta,\delta}} \sum_{k=1}^{\kappa_i} c_{i,k}\, (u_{i,k}^{(0)},p_{i,k}^{(0)}) + (V,Q),
\]
where $V\in E_\delta^2(K)$ and $Q\in V_\delta^1(K)$. Furthermore,
\begin{eqnarray} \label{4Al2} \nonumber
&& \| V\|_{E_\delta^2(k)} + \| Q\|_{V_\delta^1(K)} +  \sum_{i\in I_{\beta,\delta}} \sum_{k=1}^{\kappa_i} |c_{i,k}|  \\
&& \le c\, \Big( \| f\|_{E_\delta^0(K)} + \| g\|_{V_\delta^1(K)} + \| u\|_{E_\beta^2(K)} + \| p\|_{V_\beta^1(K)}\Big).
\end{eqnarray}
if $|s|=1$. We define
\begin{eqnarray} \label{2Al2}
(V',Q') & = & (V,Q) - \zeta_s \sum_{i\in I_{\beta,\delta}} \sum_{k=1}^{\kappa_i} \sum_{\mu=1}^{N_{i,\gamma}}
  c_{i,k} (sr^2)^\mu\, \big( u_{i,k}^{(\mu)},p_{i,k}^{(\mu)}\big)  \nonumber \\
& = & (u,p) - \zeta_s \sum_{i\in I_{\beta,\delta}} \sum_{k=1}^{\kappa_i} c_{i,k}\, \big( u_{i,k,\gamma},p_{i,k,\gamma}\big).
\end{eqnarray}
Since $\zeta_s r^{2\mu} u_{i,k}^{(\mu)}\in E_\delta^2(K)$ and $\zeta_s r^{2\mu} p_{i,k}^{(\mu)}\in V_\delta^1(K)$ for $\mu\ge 1$ we get the estimate
\begin{equation} \label{5Al2}
\| V'\|_{E_\delta^2(K)} + \| Q'\|_{V_\delta^1(K)} \le    \| V\|_{E_\delta^2(k)} + \| Q\|_{V_\delta^1(K)}
  + c\, \sum_{i\in I_{\beta,\delta}} \sum_{k=1}^{\kappa_i} |c_{i,k}|
\end{equation}
Furthermore, it follows from Lemma \ref{Al1} that
\[
(s-\Delta)V'+\nabla Q' = f -\sum_{i\in I_{\beta,\delta}} \sum_{k=1}^{\kappa_i} c_{i,k}\, \Big( (s-\Delta)\,  (\zeta_s u_{i,k,\gamma})
  + \nabla (\zeta_s p_{i,k,\gamma})\Big) \in E_\gamma^0(K)
\]
and
\[
-\nabla\cdot V' = g + \sum_{i\in I_{\beta,\delta}} \sum_{k=1}^{\kappa_i} c_{i,k}\, \nabla \cdot (\zeta_s u_{i,k,\gamma}) \in X_\gamma^1(K).
\]
Since $\gamma \ge \beta-m > \delta-m+1$, we can apply the induction hypothesis and obtain
\[
(V',Q') = \zeta_s \sum_{i\in I_{\delta,\gamma}} \sum_{k=1}^{\kappa_i} c_{i,k}\, \big( u_{i,k,\gamma},p_{i,k,\gamma}\big) + (v,q),
\]
where $(v,q)\in  E_\gamma^2(K)\times V_\gamma^1(K)$. This together with (\ref{2Al2}) implies (\ref{1Al2}). If $|s|=1$, then
it follows from the induction hypothesis that
\begin{eqnarray*}
&& \| v\|_{E_\gamma^2(k)} + \| q\|_{V_\gamma^1(K)} +  \sum_{i\in I_{\delta,\gamma}} \sum_{k=1}^{\kappa_i} |c_{i,k}|  \\
&& \le c\, \Big( \| f\|_{E_\gamma^0(K)} + \| g\|_{V_\gamma^1(K)} +   \sum_{i\in I_{\beta,\delta}} \sum_{k=1}^{\kappa_i} |c_{i,k}|
+ \| V'\|_{E_\delta^2(K)} + \| Q'\|_{V_\delta^1(K)}
  \Big).
\end{eqnarray*}
Using the last estimate together with (\ref{4Al2}) and (\ref{5Al2}), we get (\ref{3Al2}). \hfill $\Box$

\begin{Rem} \label{r1}
{\em The functions $u_{i,k}^{(\mu)}$ and $p_{i,k}^{(\mu)}$ have the form (\ref{uknu1}). We denote by $U_{i,k}^{(\mu)}(x,\tau)$ and $P_{i,k}^{(\mu)}(x,\tau)$
the functions which arise if we replace $\log r$ by $\log(\tau r)$ in (\ref{uknu1}), where $\tau$ is an arbitrary positive number, i.~e.,
\[
U_{i,k}^{(\mu)}(x,\tau) = \tau^{-\lambda_i}\, u_{i,k}^{(\mu)}(\tau x), \quad P_{i,k}^{(\mu)}(x,\tau) = \tau^{1-\lambda_i}\, p_{i,k}^{(\mu)}(\tau x).
\]
Clearly, the functions $\big( U_{i,k}^{(0)}(\cdot,\tau),P_{i,k}^{(0)}(\cdot,\tau)\big)$ form also a basis in the space of the solutions
of the problem (\ref{hom}) which have the form (\ref{sing}). Furthermore, $U_{i,k}^{(\mu)}(\cdot,\tau)$ and $P_{i,k}^{(\mu)}(\cdot,\tau)$ 
satisfy (\ref{uknu2}) and (\ref{uknu3}) for arbitrary $\tau>0$. For this reason, one can replace $u_{i,k,\gamma}(x,s)$ and $p_{i,k,\gamma}(x,s)$ in
Lemma \ref{Al2} by
\[
U_{i,k,\gamma}(x,s,\tau) = \sum_{\mu=0}^{N_{i,\gamma}} (sr^2)^\mu \, U_{i,k}^{(\mu)}(x,\tau) = \tau^{-\lambda_i}\, u_{i,k,\gamma}(\tau x,\tau^{-2}s)
\]
and
\[
P_{i,k,\gamma}(x,s,\tau) = \sum_{\mu=0}^{N_{i,\gamma}} (sr^2)^\mu \, P_{i,k}^{(\mu)}(x,\tau) = \tau^{1-\lambda_i}\, p_{i,k,\gamma}(\tau x,\tau^{-2}s).
\]
Then, instead of (\ref{1Al2}), we obtain the decomposition
\begin{equation} \label{1Al2a}
( u,p) = \zeta_s \sum_{i\in I_{\beta,\gamma}} \sum_{k=1}^{\kappa_i} C_{i,k}(s,\tau)\, \big( U_{i,k,\gamma}(x,s,\tau),P_{i,k,\gamma}(x,s,\tau)\big) + (v,q),
\end{equation}
with coefficients $C_{i,k}$ depending on $s$ and $\tau$. However, the remainder $(v,q)$ is the same as in (\ref{1Al2}).}
\end{Rem}

Indeed, the difference
\begin{equation} \label{1Al2b}
\zeta_s \sum_{i\in I_{\beta,\gamma}} \sum_{k=1}^{\kappa_i} \big( C_{i,k}(s,\tau)\, ( U_{i,k,\gamma}(x,s,\tau) - c_{i,k}(s)\, u_{i,k,\gamma}(x,s)\big)
\end{equation}
is a sum of terms of the form
\[
\zeta_s \, C_{i,k,\mu}(s,\tau)\, r^{\lambda_i+2\mu}\, a_{i,k,\mu}(\omega)\, \log^k r,
\]
where $0\le \mu\le N_{i,\gamma}$. Since $\gamma+\mbox{Re}\, \lambda_i +2N_{i,\gamma} \le \frac 12$, the expression (\ref{1Al2b}) is an element of the
space $E_\gamma^2(K)$ only if all coefficients $C_{i,k,\mu}(s,\tau)$ are zero. Hence the function $v$ in (\ref{1Al2a}) is the same as in
(\ref{1Al2}). The same holds for the function $q$.

\subsection{A formula for the coefficients in the asymptotics}

We consider the eigenvalues $\lambda_{-j}=-1-\lambda_j$ with negative real parts and define
\[
v_{j,l}^{(0)} =u_{-j,l}^{(0)}, \quad q_{j,l}^{(0)} = p_{-j,l}^{(0)} \quad \mbox{for } j=1,2,\ldots, \ l=1,2,\ldots,\kappa_j,
\]
where $\big( u_{-j,l}^{(0)},p_{-j,l}^{(0)}\big)$ are the special solutions of (\ref{hom}) introduced in the foregoing subsection.
Let
\[
a(u,p,v,q)= \int_K\Big( u\cdot (-\Delta v+ \nabla q) - p\, \nabla\cdot v\Big)\, dx
\]
and let $\chi$ be an arbitrary two times continuously differentiable function on $\overline{K}$ with compact support which is equal
to one in a neighborhood of the vertex of the cone. Then the functions  $(v_{j,l}^{(0)},q_{j,l}^{(0)})$ can be chosen such that
the biorthogonality condition
\begin{equation} \label{biorth}
a\big( u_{i,k}^{(0)},p_{i,k}^{(0)}, \chi v_{j,l}^{(0)},\chi q_{j,l}^{(0)}\big) = \delta_{i,j}\, \delta_{k,l}
\end{equation}
is satisfied for all $i,j\ge 1$, $k\le \kappa_i$, $l\le \kappa_j$ (cf. \cite[Theorem 5.1.1 and Lemma 5.1.5]{kmr-97}, where $t=\log r$).
Note that the expression on the left-hand side of (\ref{biorth}) is independent of the cut-off function $\chi$. Indeed, if $\chi'$
is another function with compact support equal to 1 near the origin, then we get
\[
a\big( u_{i,k}^{(0)},p_{i,k}^{(0)}, (\chi-\chi') v_{j,l}^{(0)},(\chi-\chi') q_{j,l}^{(0)}\big)
  = a\big( (\chi-\chi')v_{j,l}^{(0)},(\chi-\chi')q_{j,l}^{(0)}, u_{i,k}^{(0)}, p_{i,k}^{(0)}\big)=0.
\]
Let $N_j$ be the smallest integer such that $2N_j>\mbox{Re}\, \lambda_j-1$.
Analogously to the functions $u_{i,k}^{(\mu)}$, $p_{i,k}^{(\mu)}$ introduced in Subsection 2.1, let $v_{j,l}^{(\nu)}$, $q_{j,l}^{(\nu)}$,
$\nu=1,2,\ldots,N_j$, be functions of the form
\begin{equation} \label{jlnu1}
v_{j,l}^{(\nu)}(x) = r^{-\lambda_j-1} \sum_{k=0}^\sigma (\log r)^k\ a_k(\omega), \quad
q_{j,l}^{(\nu)}(x) = r^{-\lambda_j-2} \sum_{k=0}^\sigma (\log r)^k \, b_k(\omega),
\end{equation}
such that
\begin{eqnarray} \label{jlnu2}
&& - \Delta (r^{2\nu} v_{j,l}^{(\nu)}) + \nabla (r^{2\nu}q_{j,l}^{(\nu)}) = -r^{2\nu-2}v_{j,l}^{(\nu-1)}, \quad \nabla\cdot (r^{2\nu}v_{j,l}^{(\nu)})=0
  \ \mbox{ in }K, \\ \label{jlnu3} && v_{j,l}^{(\nu)}=0\ \mbox{ on }\partial K\backslash \{ 0\}.
\end{eqnarray}
The functions $v_{j,l}^{(\nu)}$, $q_{j,l}^{(\nu)}$ are uniquely determined for given $v_{j,l}^{(\nu-1)}$ if $2\nu-1-\lambda_j$
is not an eigenvalue of the pencil ${\cal L}(\lambda)$, i.~e., $\lambda_j-2\nu \not=\lambda_m$ for all $m$. If
$\lambda_j=\lambda_m+2\nu$ and $(v_{j,l}^{(\nu)},q_{j,l}^{(\nu)})$ is a solution of (\ref{jlnu2}), (\ref{jlnu3}), then every pair
\[
\big(w_{j,l}^{(\nu)}, r_{j,l}^{(\nu)}\big) = \big( v_{j,l}^{(\nu)},q_{j,l}^{(\nu)}\big)
  + r^{-2\nu}\sum_{\sigma=1}^{\kappa_m} c_\sigma \, \big( v_{m,\sigma}^{(0)},q_{m,\sigma}^{(0)} \big)
\]
is also a solution of (\ref{jlnu2}), (\ref{jlnu3}) which has the form (\ref{jlnu1}). We define
\begin{equation} \label{vqjl}
v_{j,l}(x,s) = \sum_{\nu=0}^{N_j} (sr^2)^\nu\, v_{j,l}^{(\nu)}(x), \quad q_{j,l}(x,s)= \sum_{\nu=0}^{N_j} (sr^2)^\nu\, q_{j,l}^{(\nu)}(x)
\end{equation}
for $j=1,2,\ldots$ and $l=1,\ldots,\kappa_j$. Furthermore, we introduce the bilinear form
\[
A(u,p,v,q) = a(u,p,v,q) - a(v,q,u,p).
\]

\begin{Le} \label{Al3}
Let $\chi$ be an arbitrary two times continuously differentiable function on $\overline{K}$ with compact support which is equal
to zero near the origin, and let $\gamma < \frac 12$. Furthermore, let the functions $v_{j,l}^{(0)}$ and $q_{j,l}^{(0)}$ be such that
the biorthogonality condition {\em (\ref{biorth})} is satisfied
for $0<\lambda_i,\lambda_j <\frac 12 -\gamma$, $k=1,\ldots,\kappa_i$, $l=1,\ldots,\kappa_j$. Then one can
choose functions $v_{j,l}^{(\nu)}$ and $q_{j,l}^{(\nu)}$, $\nu=1,\ldots,N_j$, of the form  {\em (\ref{jlnu1})} such that
$v_{j,l}^{(\nu)}$ and $q_{j,l}^{(\nu)}$ satisfy {\em (\ref{jlnu2})}, {\em (\ref{jlnu3})} for $\nu=1,\ldots,N_j$ and the condition
\begin{equation} \label{1Al3}
A\big(u_{i,k,\gamma},p_{i,k,\gamma},\chi v_{j,l}, \chi q_{j,l}\big) = \delta_{i,j}\, \delta_{k,l}
\end{equation}
for $0<\lambda_i,\lambda_j <\frac 12 -\gamma$, $k=1,\ldots,\kappa_i$, $l=1,\ldots,\kappa_j$.
\end{Le}

P r o o f.
Obviously,
\begin{eqnarray*}
A\big( u_{i,k,\gamma},p_{i,k,\gamma},\chi v_{j,l}, \chi q_{j,l}\big) =
  \sum_{\mu=0}^{N_{i,\gamma}} \sum_{\nu=0}^{N_j} s^{\mu+\nu} a_{i,k,\mu}^{j,l,\nu}\, ,
\end{eqnarray*}
where
\begin{equation} \label{aijkl}
a_{i,k,\mu}^{j,l,\nu} = A\big( r^{2\mu} u_{i,k}^{(\mu)},r^{2\mu}p_{i,k}^{(\mu)},\chi r^{2\nu} v_{j,l}^{(\nu)}, \chi r^{2\nu} q_{j,l}^{(\nu)}\big).
\end{equation}
It is evident that the coefficients $a_{i,k,\mu}^{j,l,\nu}$ do not depend on the cut-off function $\chi$.
Thus, one can replace $\chi$ in (\ref{aijkl}) by the function $\chi_\tau(x)=\chi({\tau}^{-1}x)$ with arbitrary positive $\tau$.
Then the substitution $x=\tau y$ in the integral representation of the form $A$ yields
\[
a_{i,k,\mu}^{j,l,\nu} = \tau^{\lambda_i-\lambda_j+2\mu+2\nu}\, a\big( r^{2\mu} U_{i,k}^{(\mu)},r^{2\mu}P_{i,k}^{(\mu)},
  r^{2\nu} V_{j,l}^{(\nu)},r^{2\nu}Q_{j,l}^{(\nu)}\big),
\]
where $U_{i,k}^{(\mu)}(x,\tau) = \tau^{-\lambda_i} \, u_{i,k}^{(\mu)}(\tau x)$, $P_{i,k}^{(\mu)}(x,\tau) = \tau^{1-\lambda_i} \, p_{i,k}^{(\mu)}(\tau x)$,
\[
V_{j,l}^{(\nu)}(x,\tau) = \tau^{\lambda_j+1} \, v_{j,l}^{(\nu)}(\tau x) \quad\mbox{and}\quad
Q_{j,l}^{(\nu)}(x,\tau) = \tau^{\lambda_j+2} \, q_{j,l}^{(\nu)}(\tau x)
\]
are polynomials of $\log\tau$. Hence, we get the representation
\[
a_{i,k,\mu}^{j,l,\nu} = \tau^{\lambda_i-\lambda_j+2\mu+2\nu}\, b_{i,k,\mu}^{j,l,\nu}(\log\tau)
\]
with a polynomial $b_{i,k,\mu}^{j,l,\nu}$. Thus,
\[
a_{i,k,\mu}^{j,l,\nu} = 0 \  \mbox{ if } \lambda_i-\lambda_j+2\mu+2\nu\not=0.
\]
Consequently,
\[
A\big( u_{i,k,\gamma},p_{i,k,\gamma},\chi v_{j,l}, \chi q_{j,l}\big) = 0 \ \mbox{ if } \lambda_j-\lambda_i \not\in\{ 0,2,4,\ldots\}.
\]
In the case $\lambda_i=\lambda_j$ (i.~e., $i=j$), we have  $a_{j,k,\mu}^{j,l,\nu}=0$ for $\mu+\nu\ge 1$ and, consequently
\begin{eqnarray*}
A\big( u_{i,k,\gamma},p_{i,k,\gamma},\chi v_{j,l}, \chi q_{j,l}\big) & = & a_{j,k,0}^{j,l,0}
  = A\big( u_{j,k}^{0)},p_{j,k}^{0}, \chi v_{j,l}^{(0)}, \chi q_{j,l}^{(0)}\big) \\
& = & a\big( u_{j,k}^{0)},p_{j,k}^{0}, \chi v_{j,l}^{(0)}, \chi q_{j,l}^{(0)}\big) = \delta_{k,l}\, .
\end{eqnarray*}
We consider the case $\lambda_j=\lambda_i+2d$, where $d$ is a positive integer. Then $N_j\ge d$, because $2N_j> \lambda_1+\mbox{Re}\, \lambda_j-1$
and $\mbox{Re}\, \lambda_i>0$. Since $a_{i,k,\mu}^{j,l,\nu}=0$ for $\mu+\nu\not=d$, we get
\begin{eqnarray*}
A\big( u_{i,k,\gamma},p_{i,k,\gamma},\chi v_{j,l}, \chi q_{j,l}\big) & = & s^d \sum_{\nu=\max(0,d-N_{i,\gamma})}^d a_{i,k,d-\nu}^{j,l,\nu}
\end{eqnarray*}
We replace the pairs $(v_{j,l}^{(d)},q_{j,l}^{(d)})$, $l=1,\ldots,\kappa_j$, by
\[
\big(w_{j,l}^{(d)}, r_{j,l}^{(d)}\big) = \big( v_{j,l}^{(d)},q_{j,l}^{(d)}\big)
  + r^{-2d}\sum_{\sigma=1}^{\kappa_i} c_{l,\sigma} \, \big( v_{i,\sigma}^{(0)},q_{i,\sigma}^{(0)} \big)
\]
(the remaining pairs $(v_{j,l}^{(\nu)},q_{j,l}^{(\nu)})$ are not changed). Then, instead of the pair $(v_{j,l},q_{j,l})$, we
obtain the pair
\[
\big( w_{j,l},r_{j,l}\big) = \sum_{\nu=0}^{N_j} (sr^2)^\nu\, \big( v_{j,l}^{(\nu)}, q_{j,l}^{(\nu)}\big)
  + s^d\, \sum_{\sigma=1}^{\kappa_i} c_{l,\sigma} \, \big( v_{i,\sigma}^{(0)},q_{i,\sigma}^{(0)} \big)
\]
Then
\begin{eqnarray*}
&& A\big( u_{i,k,\gamma},p_{i,k,\gamma},\chi w_{j,l}, \chi r_{j,l}\big) \\
&& = A\big( u_{i,k,\gamma},p_{i,k,\gamma},\chi v_{j,l}, \chi q_{j,l}\big)
  + s^d  \sum_{\sigma=1}^{\kappa_i} c_{l,\sigma} A\big( u_{i,k}^{(0)},p_{i,k}^{(0)},v_{i,\sigma}^{(0)},q_{i,\sigma}^{(0)} \big)\\
&& = A\big( u_{i,k,\gamma},p_{i,k,\gamma},\chi v_{j,l}, \chi q_{j,l}\big)   + s^d  \,c_{l,k}\, .
\end{eqnarray*}
This means, we can choose the coefficients $c_{l,\sigma}$ such that
\[
A\big( u_{i,k,\gamma},p_{i,k,\gamma},\chi w_{j,l}, \chi r_{j,l}\big) = 0 \mbox{ for }k=1,\ldots,\kappa_i,\ l=1,\ldots,\kappa_j
\]
if $\lambda_j=\lambda_i+2d$ and $d$ is an integer, $d\ge 1$. This proves the lemma. \hfill $\Box$

\begin{Rem} \label{Ar2}
{\em 1) Obviously, we can replace the cut-off function $\chi$ in Lemma \ref{Al3} by the function $\zeta_s$ introduced in the foregoing
subsection.

2) Let $U_{i,k,\gamma}(x,s,\tau)$ and $P_{i,k,\gamma}(x,s,\tau)$ be the functions introduced in Remark \ref{r1}. Furthermore, let
\[
V_{j,l}(x,s,\tau) = \tau^{\lambda_j+1}\, v_{j,l}(\tau x,\tau^{-2}s), \quad Q_{j,l}(x,s,\tau) = \tau^{\lambda_j+2}\, v_{j,l}(\tau x,\tau^{-2}s)
\]
for $\tau>0$, i.~e., $V_{j,l}$, $Q_{j,l}$ arise if we replace $\log r$ by $\log(\tau r)$ in the representation (\ref{vqjl}) of
$v_{j,l}$ and $q_{j,l}$. Then it follows from Lemma \ref{Al3} that
\[
A\big(U_{i,k,\gamma},P_{i,k,\gamma},\chi V_{j,l}, \chi Q_{j,l}\big) = \delta_{i,j}\, \delta_{k,l}.
\]
for $0<\lambda_i,\lambda_j <\frac 12 -\gamma$, $k=1,\ldots,\kappa_i$, $l=1,\ldots,\kappa_j$.}
\end{Rem}

Let $\gamma$ be a given real number less than $\frac 12$. We assume that the functions $v_{i,k}^{(\nu)}$ and $q_{i,k}^{(\nu)}$
are chosen such that the condition (\ref{1Al3}) is satisfied for $0<\lambda_i,\lambda_j <\frac 12 -\gamma$,
$k=1,\ldots,\kappa_i$, $l=1,\ldots,\kappa_j$. Analogously to Lemma \ref{Al1},
\[
(s-\Delta)\, (\zeta_s v_{i,k}) + \nabla(\zeta_s q_{i,k}) \in E_\delta^0(K) \ \mbox{ and }\
\nabla\cdot (\zeta_s v_{i,k}) \in X_\delta^1(K)
\]
for arbitrary $\delta>\mbox{Re}\, \lambda_i -2N_i-\frac 12$. Suppose that
$\max(\frac 12-\lambda_1,\mbox{Re}\, \lambda_i-2N_i-\frac 12) <\delta<\frac 12.$
Then by Theorem \ref{t1}, there exists a uniquely determined solution
$(v'_{i,k},q'_{i,k})\in E_\delta^2(K)\times X_\delta^1(K)$ of the problem
\[
(s-\Delta)\, \big(\zeta_s v_{i,k} +v'_{i,k}\big) + \nabla\big( \zeta_s q_{i,k}+ q'_{i,k}\big)=0, \ \
\nabla\cdot \big( \zeta_s v_{i,k}+v'_{i,k}\big)=0\ \mbox{ in }K,
\]
$V'_{i,k}=0$ on $\partial K\backslash\{ 0\}$ for $s\not=0$, $\mbox{Re}\, s\ge 0$. We set
\begin{equation} \label{1Al5}
v^*_{i,k}(x,s) =\zeta_s v_{i,k}(x,s) +v'_{i,k}(x,s) \quad\mbox{and}\quad q^*_{i,k}(x,s) =\zeta_s q_{i,k}(x,s) +q'_{i,k}(x,s)
\end{equation}
for $0<\lambda_i<\frac 12 - \gamma$, $k=1,\ldots,\kappa_i$.

\begin{Th} \label{At1}
Let $(u,p) \in E_\beta^2(K) \times V_\beta^1(K)$ be a solution of the problem {\em (\ref{par3})}, where
$\mbox{\em Re}\, s \ge 0$, $s\not= 0$,
\begin{equation} \label{fandg}
f\in E_\beta^0(K)\cap E_\gamma^0(K), \quad g\in X_\beta^1(K)\cap X_\gamma^1(K), \quad \gamma < \frac 12-\lambda_1 < \beta < \frac 32 + \lambda_1.
\end{equation}
We assume that the line $\mbox{\em Re}\, \lambda = \frac 12 -\gamma$ is free of eigenvalues of the pencil ${\cal L}(\lambda)$. Then
\begin{equation} \label{1At1}
(u,p) = \zeta_s \sum_{i\in I_{\gamma}} \sum_{k=1}^{\kappa_i} c_{i,k}\, \big( u_{i,k,\gamma},p_{i,k,\gamma}\big) + (v,q),
\end{equation}
where $I_{\gamma}$ is the set of all $i$ such that $0 < \mbox{\em Re}\, \lambda_i < \frac 12 - \gamma$,
\[
c_{i,k}(s) = -\int_K \big( f(x)\cdot v^*_{i,k}(x,s) + g(x)\, q^*_{i,k}(x,s)\big)\, dx
\]
are constants depending on $f$, $g$ and $s$, $(v,q)\in E_\gamma^2(K)\times V_\gamma^1(K)$ and
\begin{eqnarray} \label{5At1} \nonumber
\| v\|^2_{V_\gamma^2(K)} + |s|^2\, \| v\|^2_{V_\gamma^0(K)}+ \| q\|^2_{V_\gamma^1(K)} & \le & c\, \Big( \| f\|^2_{V_\gamma^0(K)} + \|g\|^2_{V_\gamma^1(K)}
  + |s|^{\beta-\gamma} \| f\|^2_{V_\beta^0(K)} \\
&& \ \ + \ |s|^{\beta-\gamma}\| g\|^2_{V_\beta^1(K)}  + |s|^{\beta-\gamma+2}\, \| g\|^2_{(V_{-\beta}^1(K))^*}\Big).
\end{eqnarray}
The constants $c_{i,k}$ satisfy the estimate
\begin{eqnarray} \label{4At1} \nonumber
|c_{i,k}|^2 & \le & c\, |s|^{\gamma+\mbox{\em \scriptsize Re}\, \lambda_i -1/2}\, \big( 1+\big|\log |s|\big|^{2d_i}\big) \,
  \big( \| f\|^2_{V_\gamma^0(K)} +  \| g\|^2_{V_\gamma^1(K)}\big)  \\
&&  + \ c\, |s|^{\beta+\mbox{\em \scriptsize Re}\, \lambda_i -1/2}\, \big( 1+\big|\log |s|\big|^{2d_i}\big) \,
  \big( \| f\|^2_{V_\beta^0(K)} + |s|^2\,  \| g\|^2_{(V_{-\beta}^1(K))^*}\big).
\end{eqnarray}
Here,  $c$ is independent of $f,g,s$.
\end{Th}

P r o o f.
Suppose first that $\beta<\frac 12$. Then we can choose the number $\delta$ in the interval
\[
\max(\frac 12-\lambda_1,\mbox{Re}\, \lambda_i-2N_i-\frac 12)<\delta<\frac 12
\]
such that $0\le \beta+\delta\le 2$. In this case, $v'_{i,k} \in E_\delta^2(K)\subset V_{-\beta}^0(K)$ and
\[
q'_{i,k} \in V_\delta^1(K) \subset V_{-\beta}^1(K)+V_{2-\beta}^1(K) \subset  V_{-\beta}^1(K)+\big( V_{\beta}^1(K)\big)^* = \big( X_\beta^1(K)\big)^*.
\]
Analogously, $u\in E_\beta^2(K) \subset V_{-\delta}^0(K)$ and $p\in V_\beta^1(K) \subset \big( X_\delta^1(K)\big)^*$. Thus,
\begin{eqnarray*}
&& \int_K \big( f\cdot v'_{i,k} + g\, q'_{i,k}\big)\, dx = \int_K \Big( (su-\Delta u+\nabla p)\cdot v'_{i,k} - (\nabla\cdot u)\, q'_{i,k}\Big)\, dx  \\
&& \hspace{2em} =  \int_K \Big( u\cdot \big( (s-\Delta)v'_{i,k}+\nabla q'_{i,k}\big) - p\, \nabla\cdot v'_{i,k}\Big)\, dx \\
&& \hspace{2em} = - \int_K \Big( u\cdot \big( (s-\Delta) (\zeta_s v_{i,k})+\nabla (\zeta_s q_{i,k})\big) - p\, \nabla\cdot (\zeta_s v_{i,k})\Big)\, dx.
\end{eqnarray*}
Consequently,
\begin{eqnarray*}
&& \int_K \big( f\cdot v^*_{i,k} + g\, q^*_{i,k}\big)\, dx
  = \int_K \Big( ( su-\Delta u+\nabla p\big)\cdot \zeta_s v_{i,k} - (\nabla\cdot u)\, \zeta_s q_{i,k}\Big)\, dx \\
&& \hspace{2em} - \int_K \Big( u\cdot \big( (s-\Delta) (\zeta_s v_{i,k})+\nabla (\zeta_s q_{i,k})\big) - p\, \nabla\cdot (\zeta_s v_{i,k})\Big)\, dx.
\end{eqnarray*}
Here, $(u,p)$ admits the decomposition (\ref{1At1}) with a remainder $(v,q) \in E_\gamma^2(K)\times V_\gamma^1(K)$. Since
$\zeta_s v_{i,k} \in V_{2-\gamma}^2(K)$ and $\zeta_s q_{i,\gamma}\in V_{2-\gamma}^1(K)$ for $\lambda_i < \frac 12 -\gamma$, we get
\begin{eqnarray*}
&& \int_K \Big( (sv-\Delta v+\nabla q)\cdot \zeta_s v_{i,k} - (\nabla\cdot v)\, \zeta_s q_{i,k}\Big)\, dx \\
&& \hspace{2em} - \int_K \Big( v\cdot \big( (s-\Delta) (\zeta_s v_{i,k})+\nabla (\zeta_s q_{i,k})\big) - q\, \nabla\cdot (\zeta_s v_{i,k})\Big)\, dx = 0.
\end{eqnarray*}
This implies
\begin{eqnarray*}
&& \int_K \big( f\cdot v^*_{i,k} + g\, q^*_{i,k}\big)\, dx \\
&& = \sum_{j\in I_{\gamma}} \sum_{l=1}^{\kappa_j} c_{j,l}\, \Big(
\int_K \Big( \big( -\Delta(\zeta_s u_{j,l,\gamma})+\nabla (\zeta_s p_{j,l,\gamma})\big)\cdot \zeta_s v_{i,k}
  - \big(\nabla\cdot (\zeta_s u_{j,l,\gamma})\big)\, \zeta_s q_{i,k}\Big)\, dx \\
&& \hspace{3em} - \int_K \Big( \zeta_s u_{j,l,\gamma}\cdot \big( -\Delta (\zeta_s v_{i,k})+\nabla (\zeta_s q_{i,k})\big)
  - \zeta_s p_{j,l,\gamma}\, \nabla\cdot (\zeta_s v_{i,k})\Big)\, dx \Big)\\
&& = - \sum_{j\in I_{\gamma}} \sum_{l=1}^{\kappa_j} c_{j,l} \, A\big( \zeta_s u_{j,l,\gamma},  \zeta_s p_{j,l,\gamma}, \zeta_s v_{i,k},\zeta_s q_{i,k}\big).
\end{eqnarray*}
Obviously,
\[
A\big( (1-\zeta_s) u_{j,l,\gamma}\, ,  (1-\zeta_s) p_{j,l,\gamma}\, , \zeta_s v_{i,k}\, ,\zeta_s q_{i,k}\big)=0.
\]
Therefore, Lemma \ref{Al3} implies
\[
\int_K \big( f\cdot v^*_{i,k} + g\, q^*_{i,k}\big)\, dx = -\sum_{j\in I_{\gamma}} \sum_{l=1}^{\kappa_j} c_{j,l} \,
  A\big( u_{j,l,\gamma},  p_{j,l,\gamma}, \zeta_s v_{i,k},\zeta_s q_{i,k}\big) = - c_{i,k}\, .
\]
We estimate the coefficient
\[
c_{i,k} = - \int_K \Big( f\cdot \big(\zeta_s v_{i,k}+v'_{i,k}\big) + g\, \big( \zeta_s q_{i,k}+ q'_{i,k}\big)\Big)\, dx.
\]
Obviously,
\[
\big| \zeta_s v_{i,k}(x)\big| \le c\, r^{-\mbox{\scriptsize Re}\, \lambda_i-1}\, \, \big( 1+|\log r|^{d_i}\big), \quad
\big| \zeta_s q_{i,k}(x)\big| \le c\, r^{-\mbox{\scriptsize Re}\, \lambda_i-2}\, \, \big( 1+|\log r|^{d_i}\big),
\]
where $0< \mbox{Re}\, \lambda_i < \frac 12 -\gamma$. Since moreover $\zeta_s(x)=0$ for $|sr^2|>1$, we obtain
\begin{eqnarray*}
\Big| \int_K f\, \zeta_s\, v_{i,k}\, dx\Big|^2 & \le & c \, \| f\|^2_{V_\gamma^0(K)}\ \int_0^{\sqrt{2/|s|}} r^{-2\gamma-2\, \mbox{\scriptsize Re}\, \lambda_i}\,
  \big( 1+|\log r|^{2d_i}\big)\, dr \\
& \le & c\, |s|^{\gamma+\mbox{\scriptsize Re}\, \lambda_i -1/2}\, \big( 1+\big|\log |s|\big|^{2d_i}\big)\, \| f\|^2_{V_\gamma^0(K)}
\end{eqnarray*}
and
\begin{eqnarray*}
\Big| \int_K g\, \zeta_s\, q_{i,k}\, dx\Big|^2 & \le & c \, \| g\|^2_{V_{\gamma-1}^0(K)}\ \int_0^{\sqrt{2/|s|}} r^{-2\gamma-2\, \mbox{\scriptsize Re}\, \lambda_i}\,
  \big( 1+|\log r|^{2d_i}\big)\, dr \\
& \le & c\, |s|^{\gamma+\mbox{\scriptsize Re}\, \lambda_i -1/2}\, \big( 1+\big|\log |s|\big|^{2d_i}\big)\, \| g\|^2_{V_\gamma^1(K)} \, .
\end{eqnarray*}
We consider the integral of $f\cdot v'_{i,k} + g\,  q'_{i,k}$ over the cone $K$.
The pair $\big( v'_{i,k},q'_{i,k}\big)$ is a solution of the problem
\begin{equation} \label{vqik}
(s-\Delta)\, v'_{i,k} + \nabla q'_{i,k} = F_{i,k}, \ \ -\nabla\cdot v'_{i,k} =G_{i,k}\ \mbox{ in }K, \quad v'_{i,k}=0\ \mbox{ on }\partial K\backslash\{ 0\},
\end{equation}
where
\begin{equation} \label{FGik}
F_{i,k} = -(s-\Delta)\, (\zeta_s v_{i,k}) - \nabla (\zeta_s q_{i,k}), \quad G_{i,k}= \nabla\cdot(\zeta_s v_{i,k}) = v_{i,k}\cdot \nabla\zeta_s.
\end{equation}
Furthermore, $v'_{i,k} \in E_\delta^2(K)$ and $q'_{i,k}\in V_\delta^1(K)$, where
$\max(\frac 12-\lambda_1,\mbox{Re}\, \lambda_i-2N_i-\frac 12)<\delta<\frac 12$, $\beta+\delta\ge 0$.
By Theorem \ref{t1}, the functions $v'_{i,k}$ and $q'_{i,k}$ satisfy the estimate
\begin{eqnarray} \label{2At1}
&& \| v'_{i,k}\|^2_{V_\delta^2(K)} + |s|^2\, \| v'_{i,k}\|^2_{V_\delta^0(K)} + \| q'_{i,k}\|^2_{V_\delta^1(K)} \nonumber \\
&&  \le c\, \Big( \| F_{i,k}\|^2_{V_\delta^0(K)}
  + \| G_{i,k} \|^2_{V_\delta^1(K)} + |s|^2\, \| G_{i,k}\|^2_{(V_{-\delta}^1(K))^*}\Big)
\end{eqnarray}
with a constant $c$ independent of $s$. We estimate the right-hand side of (\ref{2At1}).
Using the equality $(s-\Delta)v_{i,k}+\nabla q_{i,k} = s\, (sr^2)^{N_i}\, v_{i,k}^{(N_i)}$ (cf. (\ref{1Al1}))
and the fact that $\zeta_s(x)=0$ for $|sr^2|>2$, $\zeta_s(x)=1$ for $|sr^2|<1$, we get the estimate
\[
| F_{i,k}| \le c\, |s| \, r^{-1-\lambda_i}\, |sr^2|^{N_i}\, \big( 1+|\log r|^{d_i}\big).
\]
This implies
\begin{eqnarray*}
\| F_{i,k}\|^2_{V_\delta^0(K)} & \le & c\, |s|^{2N_i+2} \int_0^{\sqrt{2/|s|}} r^{2\delta-2\mbox{\scriptsize Re}\, \lambda_i+4N_i} \,
  \big( 1+|\log r|^{2d_i}\big)\, dr \\ & \le & c\, |s|^{-\delta+ \mbox{\scriptsize Re}\, \lambda_i +3/2} \, \big( 1+\big|\log |s|\big|^{2d_i}\big)
\end{eqnarray*}
since $2\delta-2\mbox{Re}\, \lambda_i+4N_i>-1$. The function $G_{i,k}$ satisfies the estimates
\[
\big| G_{i,k}(x)\big| \le c\, r^{-2-\lambda_i}\, \big( 1+|\log r|^{d_i}\big) \quad\mbox{and}\quad
\big| \nabla G_{i,k}(x)\big| \le c\, r^{-3-\lambda_i}\, \big( 1+|\log r|^{d_i}\big)
\]
Furthermore, $G_{i,k}(x)=0$ for $|sr^2|<1$ and $sr^2>2$. Thus, one easily obtains the estimates
\[
\| G_{i,k} \|^2_{V_\delta^1(K)} \le c\, |s|^{-\delta+ \mbox{\scriptsize Re}\, \lambda_i +3/2} \, \big( 1+\big|\log |s|\big|^{2d_i}\big)
\]
and
\[
\| G_{i,k}\|^2_{(V_{-\delta}^1(K))^*} \le \| G_{i,k}\|^2_{V_{\delta+1}^0(K)}
  \le c\, |s|^{-\delta+ \mbox{\scriptsize Re}\, \lambda_i -1/2} \, \big( 1+\big|\log |s|\big|^{2d_i}\big).
\]
Thus, by (\ref{2At1}),
\begin{equation} \label{3At1}
\| v'_{i,k}\|^2_{V_\delta^2(K)} + |s|^2\, \| v'_{i,k}\|^2_{V_\delta^0(K)} + \| q'_{i,k}\|^2_{V_\delta^1(K)}
  \le c\, |s|^{-\delta+ \mbox{\scriptsize Re}\, \lambda_i +3/2} \, \big( 1+\big|\log |s|\big|^{2d_i}\big).
\end{equation}
Let $K_s= \{ x\in K:\, |sr^2|<1\}$. Since $\beta+\delta\ge 0$ and $\gamma+\delta\le 2$, we obtain
\begin{eqnarray*} 
&& \hspace{-1em} \frac 12\, \Big| \int_K f\cdot v'_{i,k} \, dx\Big|^2
\le  \| f\|^2_{V_\gamma^0(K)}\, \int_{K_s} r^{-2\gamma}\, |v'_{i,k}|^2\, dx
  + \| f\|^2_{V_\beta^0(K)}\, \int_{K\backslash K_s} r^{-2\beta}\, |v'_{i,k}|^2\, dx \\ \nonumber
&& \le   |s|^{\gamma+\delta-2}\, \| f\|^2_{V_\gamma^0(K)}\,\int_{K_s} r^{2\delta-4}\, |v'_{i,k}|^2\, dx
  + |s|^{\beta+\delta}\, \| f\|^2_{V_\beta^0(K)}\,\int_{K\backslash K_s} r^{2\delta}\, |v'_{i,k}|^2\, dx \\
&& \le  |s|^{\delta-2}\, \Big( |s|^{\gamma}\, \| f\|^2_{V_\gamma^0(K)} + |s|^{\beta}\, \| f\|^2_{V_\beta^0(K)}\Big) \
  \Big( \| v'_{i,k}\|^2_{V_\delta^2(K)} + |s|^2\, \| v'_{i,k}\|^2_{V_\delta^0(K)}\Big).
\end{eqnarray*}
Analogously,
\begin{eqnarray*} 
&& \hspace{-1em} \frac 12 \,\Big| \int_K g\, q'_{i,k} \, dx\Big|^2
\le    \| g\|^2_{V_{\gamma-1}^0(K)}\, \int_{K_s} r^{2-2\gamma}\, |q'_{i,k}|^2\, dx
  +  \| g\|^2_{V_{\beta+1}^0(K)}\, \int_{K\backslash K_s} r^{-2\beta-2}\, |q'_{i,k}|^2\, dx \\ \nonumber
&& \le   |s|^{\gamma+\delta-2}\, \| g\|^2_{V_{\gamma-1}^0(K)}\,\int_{K_s} r^{2\delta-2}\, |v'_{i,k}|^2\, dx
  + |s|^{\beta+\delta}\, \| g\|^2_{V_{\beta+1}^0(K)}\,\int_{K\backslash K_s} r^{2\delta-2}\, |v'_{i,k}|^2\, dx \\
&& \le   |s|^{\delta-2}\, \Big( |s|^{\gamma}\, \| g\|^2_{V_\gamma^1(K)} + |s|^{\beta+2}\, \| g\|^2_{(V_{-\beta}^1(K))^*}\Big) \,
   \| q'_{i,k}\|_{V_\delta^1(K)}
\end{eqnarray*}
Using (\ref{3At1}), we obtain
\begin{eqnarray*}
&& \hspace{-1em} \Big| \int_K \big( f\cdot v'_{i,k} +g\, q'_{i,k}\big)\, dx\Big|^2
  \le   c\, |s|^{\gamma+\mbox{\scriptsize Re}\, \lambda_i -1/2}\, \big( 1+\big|\log |s|\big|^{2d_i}\big) \,
  \big( \| f\|^2_{V_\gamma^0(K)} +  \| g\|^2_{V_\gamma^1(K)}\big) \\
&& \qquad + \  c\, |s|^{\beta+\mbox{\scriptsize Re}\, \lambda_i -1/2}\, \big( 1+\big|\log |s|\big|^{2d_i}\big) \,
  \big( \| f\|^2_{V_\beta^0(K)} + |s|^2\,  \| g\|^2_{(V_{-\beta}^1(K))^*}\big) 
\end{eqnarray*}
This proves (\ref{4At1}). We estimate the norms of $v$ and $q$. In the case $|s|=1$, the estimate (\ref{5At1})
follows directly from Theorem \ref{t1} and Lemma \ref{Al2}. If $|s|$ is arbitrary, we set $x=|s|^{-1/2} y$ and define
\[
\hat{u}(y)=u(x), \ \ \hat{p}(y)=|s|^{-1/2}\, p(x), \ \ \hat{f}(y)=|s|^{-1}\, f(x),\ \ \hat{g}(y)=|s|^{-1/2}\, g(x).
\]
Obviously, $(\hat{u},\hat{p})$ is a solution of the Dirichlet problem for the system
\begin{equation} \label{9At1}
\big( |s|^{-1}s - \Delta\big)\, \hat{u} + \nabla\, \hat{p}=\hat{f}, \ \ -\nabla\cdot \hat{u}=\hat{g}\ \mbox{ in }K.
\end{equation}
Consequently,
\[
\big( \hat{u}(y),\hat{p}(y)\big) = \chi(|y|^2) \sum_{i\in I_{\gamma}} \sum_{k=1}^{\kappa_i} c_{i,k}(|s|^{-1}s)\,
  \big( u_{i,k,\gamma}(y,|s|^{-1}s),p_{i,k,\gamma}(y,|s|^{-1}s)\big)   + \big(\hat{v}(y),\hat{q}(y)\big),
\]
where
\begin{equation} \label{6At1}
\| \hat{v}\|^2_{E_\gamma^2(k)} + \| \hat{q}\|^2_{V_\gamma^1(K)}
\le c\, \Big( \| \hat{f}\|^2_{E_\gamma^0(K)} + \| \hat{g}\|^2_{V_\gamma^1(K)} + \| \hat{f}\|^2_{E_\beta^0(K)} + \| \hat{g}\|^2_{X_\beta^1(K)}\Big).
\end{equation}
By Remark \ref{r1}, we can replace $u_{i,k,\gamma}(y,|s|^{-1}s)$, $p_{i,k,\gamma}(y,|s|^{-1}s)$ by the functions $U_{i,k,\gamma}(y,|s|^{-1}s,\tau)$ and
$P_{i,k,\gamma}(y,|s|^{-1}s,\tau)$, respectively, and obtain the decomposition
\[
\big( \hat{u}(y),\hat{p}(y)\big) = \zeta(|y|^2) \sum_{i\in I_{\gamma}} \sum_{k=1}^{\kappa_i} C_{i,k}(|s|^{-1}s,\tau)\,
  \big( U_{i,k,\gamma}(y,|s|^{-1}s,\tau),P_{i,k,\gamma}(y,|s|^{-1}s,\tau)\big)   + \big(\hat{v}(y),\hat{q}(y)\big)
\]
with the same remainder $\big(\hat{v}(y),\hat{q}(y)\big)$.
Substituting $y=|s|^{1/2}x$ and $\tau=|s|^{-1/2}$, we get
\[
\big( u(x),p(x)\big) = \zeta_s(x) \sum_{k=1}^{\kappa_i} |s|^{\lambda_i/2}\, C_{i,k}(|s|^{-1}s,|s|^{-1/2})\,
\big( u_{i,k,\gamma}(x,s),p_{i,k,\gamma}(x,s)\big)  + \big(v(x),q(x)\big),
\]
where $v(x)=\hat{v}(y)$ and $q(x)=|s|^{1/2}\hat{q}(y)$. Using (\ref{6At1}) and the equalities
\[
\| v\|^2_{V_\gamma^2(K)} + |s|^2 \, \| v\|^2_{V_\gamma^0(K)} + \| q\|^2_{V_\gamma^1(K)} =
  |s|^{-\gamma+1/2}\, \Big( \| \hat{v}\|^2_{E_\gamma^2(K)} + \| \hat{q}\|^2_{V_\gamma^1(K)}\Big)
\]
and
\[
\| f\|^2_{V_\beta^0(K)} + \| g\|^2_{V_\beta^1(K)} + |s|^2\, \| g\|^2_{(V_{-\beta}^1(K))^*}
  = |s|^{-\beta+1/2} \, \Big( \| \hat{f}\|^2_{V_\beta^0(K)} + \|\hat{g}\|^2_{X_\beta^1(K)} \Big),
\]
we get (\ref{5At1}). This proves the lemma for the case $\beta<\frac 12$.

Suppose now that $\frac 12 \le \beta <\frac 32 + \lambda_1$ and that $\beta'$ is a number in the interval $(\frac 12 -\lambda_1\, ,\, \frac 12 )$.
Obviously, $f\in E_{\beta'}^0(K)$ and $g\in X_{\beta'}^1(K)$. Since $\lambda_{-1}= -1-\lambda_1 < \frac 12 -\beta < \frac 12 -\beta' < \lambda_1$,
it follows from Lemma \ref{l5b} that $u \in E_{\beta'}^2(K)$ and $p \in V_{\beta'}^1(K)$. Consequently, by the first part of the proof, we get the representation
(\ref{1At1}), where $(v,q)$ and $c_{i,k}$ satisfy the estimates (\ref{5At1}) and (\ref{4At1}) with $\beta'$ instead of $\beta$.
Using the inequality $r^{2\beta'} < |s|^{\beta-\beta'}\, r^{2\beta} + |s|^{\gamma-\beta'}\, r^{2\gamma}$ for $\gamma<\beta'<\beta$, we obtain the
estimate
\[
\| f\|^2_{V_{\beta'}^0(K)} \le |s|^{\beta-\beta'}\, \| f\|^2_{V_{\beta}^0(K)} + |s|^{\gamma-\beta'}\, \| f\|^2_{V_{\gamma}^0(K)} \, .
\]
Analogously,
\[
\| g\|^2_{V_{\beta'}^1(K)} \le |s|^{\beta-\beta'}\, \| g\|^2_{V_{\beta}^1(K)} + |s|^{\gamma-\beta'}\, \| g\|^2_{V_{\gamma}^1(K)}
\]
and
\[
\| g\|^2_{(V_{-\beta'}^1(K))^*} \le |s|^{\beta-\beta'}\, \| g\|^2_{(V_{-\beta}^1(K))^*} + |s|^{\gamma-2-\beta'}\, \| g\|^2_{V_{\gamma}^1(K)}\, .
\]
Thus, $(v,q)$ and the coefficients $c_{i,k}$ satisfy the estimates (\ref{5At1}) and (\ref{4At1}), respectively.
The proof of the theorem is complete. \hfill $\Box$ \\

Note that the number $d_i$ in the estimate (\ref{4At1}) is the greatest exponent of $\log r$ in the representation of $(v_{i,k},q_{i,k})$.

\begin{Rem} \label{Ar3}
{\em Under the conditions of Theorem \ref{At1}, the solution $(u,p)$ in the last theorem has also the representation
\[
( u,p) = \zeta_s \sum_{i\in I_\gamma} \sum_{k=1}^{\kappa_i} C_{i,k}(s,\tau)\, \big( U_{i,k,\gamma}(x,s,\tau),P_{i,k,\gamma}(x,s,\tau)\big) + (v,q)
\]
(cf. (\ref{1Al2a})) for arbitrary $\tau>0$, where $U_{i,k,\gamma}(x,s,\tau)$, $P_{i,k,\gamma}(x,s,\tau)$ are the functions introduced in Remark \ref{r1}
and $(v,q)$ is the same remainder as in (\ref{1At1}). The coefficients $C_{i,k}$ are given by the formula
\[
C_{i,k}(s,\tau) = -\int_K \big( f(x)\cdot V^*_{i,k}(x,s,\tau) + g(x)\, Q^*_{i,k}(x,s,\tau)\big)\, dx,
\]
where $V^*_{i,k}(x,s,\tau)= \tau^{\lambda_i+1}\, v^*_{i,k}(\tau x,\tau^{-2}s)$ and $Q^*_{i,k}(x,s,\tau)= \tau^{\lambda_i+2}\, q^*_{i,k}(\tau x,\tau^{-2}s)$
(cf. Remark \ref{Ar2}).}
\end{Rem}

Our goal is to estimate the coefficients $c_{i,k}$ and the remainder $(v,q)$ in (\ref{1At1}) only by the $V_\gamma^0$-norm of $f$ and
the $X_\gamma^1$-norm of $g$. For this, we need estimates of the functions $v^*_{i,k}$ and $q^*_{i,k}$ which appear in the formula
for the coefficients $c_{i,k}$. We introduce the following weighted H\"older space $N_\beta^{l,\sigma}(K)$ with the norm
\[
\| u\|_{N_\beta^{l,\sigma}(K)} = \sum_{|\alpha|\le l} \sup_{x\in K} |x|^{\beta-l-\sigma+|\alpha|}\, \big| \partial_x^\alpha u(x)\big|
  + \sum_{|\alpha|=l} \sup_{\substack{x,y\in K \\ 2|x-y|<|x|}} |x|^\beta \, \frac{|\partial_x^\alpha u(x) - \partial_y^\alpha u(y)}{|x-y|^\sigma}\, ,
\]
where $l$ is a nonnegative integer, $\beta$ and $\sigma$ are real numbers, $0<\sigma<1$. Note that
\[
V_\beta^{l+2}(K) \subset N_{\beta+\sigma-1/2}^{l,\sigma}(K) \quad\mbox{for }\sigma<1/2
\]
(cf. \cite[Lemma 3.6.2]{mr-10}).

\begin{Le} \label{Al4}
Suppose that $(u,p)\in E_\delta^2(K)\times V_\delta^1(K)$ is a solution of the problem {\em (\ref{par3})}, where
\[
f\in V_\delta^0(K) \cap N_{\delta+\sigma+3/2}^{0,\sigma}(K), \ \ g\in V_\delta^1(K) \cap N_{\delta+\sigma+3/2}^{1,\sigma}(K), \quad
  \frac 12- \lambda_1 < \delta < \frac 12 \, , \ 0<\sigma <\frac 12\, .
\]
If $f(x)$ and $g(x)$ are zero for $|s|\, r^2>1$, then $u(x)$ and $p(x)$ satisfy the estimate
\begin{equation} \label{1Al4}
\sum_{|\alpha|\le 2} r^{\delta+|\alpha|-1/2}\, \big| \partial_x^\alpha u(x)\big| + \sum_{|\alpha|\le 1} r^{\delta+|\alpha|+1/2}\,
  \big| \partial_x^\alpha p(x)\big|   \le c\, \| (f,g)\|_\delta
\end{equation}
for $|s|\, r^2 <1$, where
\[
\| (f,g)\|_\delta = \| f\|_{V_\delta^0(K)} + \| f\|_{N_{\delta+\sigma+3/2}^{0,\sigma}(K)} + \| g\|_{V_\delta^1(K)} + \| g\|_{N_{\delta+\sigma+3/2}^{1,\sigma}(K)}\, .
\]
Furthermore, the estimate
\[
\sum_{|\alpha|\le 2} |s|^{(\mu_2+2-|\alpha|)/2}\, r^{\mu_2+2}\, \big| \partial_x^\alpha u(x)\big| + \sum_{|\alpha|\le 1} |s|^{\mu_2/2}\, r^{\mu_2+1+|\alpha|}\,
  \big| \partial_x^\alpha p(x)\big|   \le c\, |s|^{(\delta-1/2)/2}\, \| (f,g)\|_\delta
\]
holds for $|s|\, r^2 >1$ if $g$ satisfies the condition {\em (\ref{2t1})}. If {\em (\ref{2t1})} is not satisfied,
then the eigenvalue $\mu_2$ in the last estimate has to be replaced by $\mu_1=0$.
\end{Le}

P r o o f. We assume first that $|s|=1$. Let $\chi$ be an infinitely differentiable function in $K$, $\chi(x)=1$ for $|x|<1$, $\chi(x)=0$ for $|x|>2$.
Then
\[
-\Delta (\chi u)+ \nabla(\chi p) = f - [\Delta,\chi] u + p\, \nabla \chi - s\chi u = F, \quad -\nabla\cdot (\chi u) = g - u\cdot\nabla\chi = G
\]
in $K$, where $[\Delta,\chi]\, u = \Delta(\chi u)-\chi\Delta u$. Here, $\chi u \in V_{\delta+2}^2(K) \subset N_{\delta+\sigma+3/2}^{0,\sigma}(K)$.
Furthermore, it follows from well-known regularity results for elliptic systems (see \cite{adn}) that
$[\Delta,\chi] u + p\, \nabla \chi \in N_{\delta+\sigma+3/2}^{0,\sigma}(K)$, $u\cdot\nabla\chi \in N_{\delta+\sigma+3/2}^{1,\sigma}(K)$ and
\begin{eqnarray*}
&& \big\| [\Delta,\chi] u + p\, \nabla \chi \big\|_{N_{\delta+\sigma+3/2}^{0,\sigma}(K)}
  + \big\| u\cdot\nabla\chi \|_{N_{\delta+\sigma+3/2}^{1,\sigma}(K)}  \\ && \le c\, \Big( \| f\|_{N_{\delta+\sigma+3/2}^{0,\sigma}(K)}
   + \| g\|_{N_{\delta+\sigma+3/2}^{1,\sigma}(K)} + \| u\|_{V_\delta^2(K)} + \| p\|_{V_\delta^1(K)}\Big).
\end{eqnarray*}
Applying \cite[Theorem 5.1, Corollary 5.1]{mp-78a}, we obtain the estimate
\[
\| \chi u\|_{N_{\delta+\sigma+3/2}^{2,\sigma}(K)} + \| \chi p\|_{N_{\delta+\sigma+3/2}^{2,\sigma}(K)}
  \le c\, \Big( \| F\|_{N_{\delta+\sigma+3/2}^{0,\sigma}(K)} + \| G\|_{N_{\delta+\sigma+3/2}^{2,\sigma}(K)}\Big).
\]
Here, by Theorem \ref{t1},
\begin{eqnarray*}
&& \| F\|_{N_{\delta+\sigma+3/2}^{0,\sigma}(K)} + \| G\|_{N_{\delta+\sigma+3/2}^{2,\sigma}(K)}  \\
&& \le  c\, \Big( \| f\|_{N_{\delta+\sigma+3/2}^{0,\sigma}(K)} + \| g\|_{N_{\delta+\sigma+3/2}^{2,\sigma}(K)} + \| u\|_{V_\delta^2(K)} + \| p\|_{V_\delta^1(K)}\Big)
 \le  c\, \| (f,g)\|_\delta\, .
\end{eqnarray*}
In particular, the estimate (\ref{1Al4}) holds for $|x|<1$, where $c$ is independent of $s$.

We estimate $u(x)$ and $p(x)$ for $|x|>1$, $|s|=1$. Since $f(x)$ and $g(x)$ are zero for $|x|>1$, we have
$f\in E_\gamma^0(K)\cap N_{\gamma+\sigma+3/2}^{0,\sigma}(K)$ and $g \in X_\gamma^0(K)\cap N_{\gamma+\sigma+3/2}^{1,\sigma}(K)$ for arbitrary $\gamma>\delta$.
Let $\gamma$ be such that $\gamma > \delta+2$ and $\gamma-\frac 12$ is not an eigenvalue of the pencil ${\cal N}(\lambda)$. We
denote the set of all integer $j$ such that $0\le \mu_j<\gamma-\frac 12$ by $J_\gamma$. By \cite[Lemma 2.5]{k/r-18a},
the vector function $(u,p)$ admits the decomposition
\begin{equation} \label{decup}
(u,p) = \eta \sum_{j\in J_\gamma} \sum_{l=1}^{\sigma_j}  c_{j,l} \, \big( u^{j,l,\gamma},p^{j,l,\gamma}\big) + (v,q).
\end{equation}
Here, $\eta=1-\zeta$, the functions $u^{j,l,\gamma}$ and $p^{j,l,\gamma}$ satisfy the inequalities
\begin{equation} \label{3Al4}
|\eta\partial_x^\alpha u^{j,l,\gamma}(x)| \le c\, r^{-2-\mu_j} \ \mbox{ for }|\alpha|\le 2, \quad
|\eta\partial_x^\alpha p^{j,l,\gamma}(x)| \le c\, r^{-1-\mu_j-|\alpha|} \ \mbox{ for }|\alpha|\le 1,
\end{equation}
and the remainder $(v,q)$ satisfies the estimate
\[
\| v\|_{E_\gamma^2(K)} + \| q\|_{V_\gamma^1(K)} + \sum_{j\in J_\gamma} \sum_{l=1}^{\sigma_j} |c_{j,l}| \le
    c\, \Big(  \| f\|_{E_\gamma^0(K)} + \| g\|_{V_\gamma^1(K)} + \| u\|_{E_\delta^2(K)} + \| p\|_{V_\delta^1(K)}\Big),
\]
where $c$ is independent of $s$. If $g$ satisfies the condition (\ref{2t1}), then it follows from \cite[Lemma 1.13]{k/r-16} that
$u\in E_{\delta'}^2(K)$ and $p\in V_{\delta'}^1(K)$ with arbitrary $\delta'$, $\frac 12 <\delta'<\mu_2+\frac 12$.
In this case, the set $J_\gamma$ in (\ref{decup}) can be replaced by $J_\gamma \backslash \{ 1\}$.
Theorem \ref{t1} implies
\[
\| v\|_{E_\gamma^2(K)} + \| q\|_{V_\gamma^1(K)} + \sum_{j\in J_\gamma} \sum_{l=1}^{\sigma_j} |c_{j,l}| \le c\, \| (f,g) \|_\delta\, .
\]
Since $(v,q)=(u,p)$ in the neighborhood $|x|<1/2$ of the origin, we have $(v,q)\in V_{\gamma-2}^2(K)\times
V_{\gamma-2}^1(K)$ if $\gamma\ge \delta+2$. Furthermore,
\[
-\Delta\, v + \nabla q = f-f'-sv, \quad -\nabla \cdot v = g-g' \ \mbox{ in }K, \quad v=0\ \mbox{ on }\partial K\backslash\{ 0\},
\]
where
\[
f' = \sum_{j\in J_\gamma} \sum_{l=1}^{\sigma_j} c_{j,l} \Big( (s-\Delta)\, ( \eta  \,  u^{j,l,\gamma})   + \nabla ( \eta p^{j,l,\gamma})\Big),\quad
g' = -\sum_{j\in J_\gamma} \sum_{l=1}^{\sigma_j} c_{j,l} \nabla\cdot ( \eta  \,  u^{j,l,\gamma}).
\]
It follows from \cite[Lemmas 2.4 and 2.5]{k/r-18a} that $f'\in V_{\gamma-2}^0(K) \cap N_{\gamma+\sigma-1/2}^{0,\sigma}(K)$,
$g'\in V_{\gamma-2}^1(K) \cap N_{\gamma+\sigma-1/2}^{1,\sigma}(K)$ and
\[
\| f'\|_{N_{\gamma+\sigma-1/2}^{0,\sigma}(K)} + \| g'\|_{N_{\gamma+\sigma-1/2}^{1,\sigma}(K)} \le c\, \sum_{j\in J_\gamma} \sum_{l=1}^{\sigma_j} |c_{j,l}|
  \le c'\, \| (f,g)\|_\delta\, .
\]
Applying \cite[Theorem 5.1, Corollary 5.1]{mp-78a} and the inequality
\[
\| v\|_{N_{\gamma+\sigma-1/2}^{0,\sigma}(K)}\le c\,\| v\|_{V_\gamma^2(K)}\, ,
\]
we obtain
\begin{eqnarray*}
&& \| v\|_{N_{\gamma+\sigma-1/2}^{2,\sigma}(K)} + \| q \|_{N_{\gamma+\sigma-1/2}^{1,\sigma}(K)}  \\
&& \le c\, \Big( \| f-f'-sv\|_{N_{\gamma+\sigma-1/2}^{2,\sigma}(K)} + \| g-g' \|_{N_{\gamma+\sigma-1/2}^{1,\sigma}(K)}\Big) \le c\, \| (f,g)\|_\delta
\end{eqnarray*}
with a constant $c$ independent of $s$. In particular,
\[
\sum_{|\alpha|\le 2} r^{\gamma+|\alpha|-5/2}\, \big| \partial_x^\alpha v(x)\big| + \sum_{|\alpha|\le 1} r^{\gamma+|\alpha|-3/2}\,
  \big| \partial_x^\alpha p(x)\big| \le c\, \| (f,g)\|_\delta
\]
for all $x\in K$. We may assume that $\gamma\ge \mu_2+\frac 92$. Then the last estimate together with (\ref{3Al4}) implies
\[
\sum_{|\alpha|\le 2}  r^{\mu_2+2}\, \big| \partial_x^\alpha u(x)\big| + \sum_{|\alpha|\le 1}  r^{\mu_2+1+|\alpha|}\, \big| \partial_x^\alpha p(x)\big|
  \le c\, \| (f,g)\|_\delta \ \mbox{ for }|x|>1
\]
if $g$ satisfies (\ref{2t1}). The same estimate with $\mu_1=0$ instead of $\mu_2$ holds if (\ref{2t1}) is not satisfied. Thus, the lemma is proved for $|s|=1.$

If $|s|$ is arbitrary, we consider the same functions $\hat{u},\hat{p},\hat{f},\hat{g}$ as in the proof of Theorem \ref{At1}.
The pair $(\hat{u},\hat{p})$ is a solution of the system (\ref{9At1}), where $\hat{f}(x)=0$ and $\hat{g}(x)=0$ for $|x|>1$.
Using the equalities
\begin{eqnarray*}
\big|\partial_x^\alpha u(x)\big| = |s|^{|\alpha|/2} \, \big|\partial_y^\alpha \hat{u}(y)\big|, \quad
|\partial_x^\alpha p(x)\big| = |s|^{(|\alpha|+1)/2} \, \big|\partial_y^\alpha \hat{p}(y)\big|
\end{eqnarray*}
for $y=|s|^{1/2}x$ and $\| (f,g)\|_\delta = |s|^{(-\delta+1/2)/2}\, \| (\hat{f},\hat{g}\|_\delta$, we obtain the estimate (\ref{1Al4})
for $|s|\,r^2<1$. Analogously, the desired estimate for the case $|s|\,r^2>1$ holds. \hfill $\Box$ \\

We estimate the functions (\ref{1Al5}) by means of the last lemma.

\begin{Le} \label{Al5}
Let $i\ge 1$, $k\le \kappa_i$, and let $(v_{i,k}^*,q^*_{i,k})$ be the solution of the problem {\em (\ref{hom})} which has the form {\em (\ref{1Al5})} with
a remainder $(v'_{i,k},q'_{i,k}) \in E_\delta^2(K)\times V_\delta^1(K)$, $\frac 12 - \lambda_1 < \delta < \frac 12$. If $\lambda_i\not=2m+1$
with an integer $m\ge 0$, then
\begin{eqnarray*}
\big| \partial_x^\alpha v'_{i,k}(x)\big| \le c\, |s|^{(\mbox{\em \scriptsize Re}\,\lambda_i+|\alpha|+1)/2} \, \big( 1+|\log |s||\big)^{d_i} \,
  \Big( \frac{|s|^{1/2} r}{1+|s|^{1/2}r}\Big)^{-\delta-|\alpha|+1/2}\,  \big( 1+|s|^{1/2}r\big)^{-2-\mu_2}
\end{eqnarray*}
for $|\alpha|\le 2$ and
\begin{eqnarray*}
\big| \partial_x^\alpha q'_{i,k}(x)\big| \le c\, |s|^{(\mbox{\em \scriptsize Re}\,\lambda_i+|\alpha|+2)/2} \, \big( 1+|\log |s||\big)^{d_i} \,
  \Big( \frac{|s|^{1/2} r}{1+|s|^{1/2}r}\Big)^{-\delta-|\alpha|-1/2}\,  \big( 1+|s|^{1/2}r\big)^{-1-\mu_2-|\alpha|}
\end{eqnarray*}
for $|\alpha|\le 1$. Here, $d_i$ is the greatest exponent of $\log r$ in the representations of $v_{i,k}$ and $q_{i,k}$.
If $\lambda_i=2m+1$, then the same estimates with $\mu_1=0$ instead of $\mu_2$ are valid.
\end{Le}

P r o o f.
The pair $(v'_{i,k},q'_{i,k})$ is the uniquely determined solution of the problem (\ref{vqik}), with the right-hand sides
$F_{i,k}$ and $G_{i,k}$ defined by (\ref{FGik}). Let $\|\cdot\|_\delta$ be the norm introduced in Lemma \ref{Al4}. Since
\[
(s-\Delta)\, v_{i,k} + \Delta q_{i,k} = s\, (sr^2)^{N_i}\, v_{i,k}^{(N_i)}, \quad \nabla\cdot v_{i,k} =0 \ \mbox{ in }K
\]
(cf. (\ref{1Al1})), it follows that
\[
\| (F_{i,k},G_{i,k})\|_\delta  \le c\, |s|^{(\mbox{\scriptsize Re}\, \lambda_i -\delta + 3/2)/2}\, \big( 1+|\log |s||\big)^{di}
\]
with a constant $c$ independent of $s$. We show that
\begin{equation} \label{2Al5}
\int_K G_{i,k}(x)\, dx = \int_K \nabla\cdot (\zeta_s v_{i,k})\, dx =0
\end{equation}
for $\lambda_i\not=2m+1$. Let $\tau$ be  positive and $\zeta_\tau(x)=\zeta(\tau r^2)$. Since $\nabla\cdot(\zeta_\tau v_{i,k}) = v_{i,k}\cdot\nabla \zeta_\tau$
vanishes in the regions $\tau|x|^2< 1/2$ and $\tau|x|^2 >1$, we get
\[
\int_K \nabla\cdot (\zeta_s v_{i,k})\, dx = \int\limits_{\substack{ K \\ 1/2 < \tau r^2<1}} \nabla\cdot (\zeta_\tau v_{i,k})\, dx
  = \int\limits_{K\cap S_{1/\sqrt{\tau}}} v_{i,k}\cdot\vec{n}\, d\sigma
\]
for arbitrary $\tau>0$, where $S_{1/\sqrt{\tau}}$ is the sphere with radius $\tau^{-1/2}$ and $\vec{n}$ is the normal vector to this sphere.
The vector function $v_{i,k}$ has the form
\[
v_{i,k}(x) = \sum_{\nu=0}^{N_i} (sr^2)^\nu \, v_{i,k}^{(\nu)}(x) = \sum_{\nu=0}^{N_i} (sr^2)^\nu \, r^{-\lambda_i-1}\sum_{j=0}^{\sigma_{i,k,\nu}}
  (\log r)^j\, \phi_{i,k,\nu,j}(\omega).
\]
Hence
\begin{eqnarray*}
\int_K G_{i,k}(x)\, dx
= \sum_{\nu=0}^{N_i} \sum_{j=0}^{\sigma_{i,k,\nu}} s^\nu\,
  \tau^{(\lambda_i-2\nu-1)/2} \, (-\log\sqrt{\tau})^j \int_\Omega \phi_{i,k,\nu,j}(\omega) \cdot\vec{n}\, d\omega.
\end{eqnarray*}
If $\lambda_i \not=2\nu+1$ for $\nu=0,1,\ldots,N_i$, then the last expression is either zero or a nonconstant function of $\tau$.
Since this expression must be independent of $\tau$, we obtain (\ref{2Al5}) for $\lambda_i\not= 2m+1$. Now the estimates for $v'_{i,k}$
and $q'_{i,k}$ follow directly from Lemma \ref{Al4}. \hfill $\Box$

\begin{Rem}
{\em By \cite[Theorem 5.3.2]{kmr-01}, there are no generalized eigenfunctions corresponding to the eigenvalues
of the pencil ${\cal L}(\lambda)$ in the strip $-2<\mbox{Re}\, \lambda < 1$. Furthermore,
$N_i=0$ for $\lambda_i<1$. Consequently, $v_{i,k}=v_{i,k}^{(0)}$ and $q_{i,k}=q_{i,k}^{(0)}$ do not contain logarithmic terms,
i.~e., $d_i=0$ if $0<\lambda_i<1$.

For $\lambda_i=1$, we have $N_i=1$. If there are no generalized eigenvectors corresponding to this eigenvalue, then the corresponding
functions $v_{i,k}^{(0)}$ and $q^{(0)}_{i,k}$ do not contain logarithmic terms. The same is true then for the functions
$v_{i,k}^{(1)}$ and $q_{i,k}^{(1)}$ since $\lambda_{-i}+2=1-\lambda_i=0$ is not an eigenvalue of the pencil ${\cal L}(\lambda)$.
This means that $d_i=0$ for $\lambda_i=1$ if there are no generalized eigenvectors corresponding to this eigenvalue.}
\end{Rem}

Using the representation (\ref{jlnu1}) for the functions $v_{i,k}^{(\nu)}$ and $q_{i,k}^{(\nu)}$, we easily deduce the following estimates
from Lemma \ref{Al5}.

\begin{Co} \label{Ac1}
Let $(v_{i,k}^*,q^*_{i,k})$ be the same vector function as in Lemma {\em \ref{Al5}}. If $\lambda_i\not=2m+1$
for integer $m\ge 0$, then
\begin{eqnarray*}
\big| \partial_x^\alpha v^*_{i,k}(x)\big| \le c\, |s|^{(\mbox{\em \scriptsize Re}\,\lambda_i+|\alpha|+1)/2} \, \big( 1+|\log r|\big)^{d_i} \,
  \Big( \frac{|s|^{1/2} r}{1+|s|^{1/2}r}\Big)^{-\mbox{\em \scriptsize Re}\,\lambda_i-|\alpha|-1}\,  \big( 1+|s|^{1/2}r\big)^{-2-\mu_2+\varepsilon}
\end{eqnarray*}
for $|\alpha|\le 2$ and
\[
\big| \partial_x^\alpha q^*_{i,k}(x)\big| \le c\, |s|^{(\mbox{\em \scriptsize Re}\,\lambda_i+|\alpha|+2)/2} \, \big( 1+|\log r|\big)^{d_i} \,
  \Big( \frac{|s|^{1/2} r}{1+|s|^{1/2}r}\Big)^{-\mbox{\em \scriptsize Re}\,\lambda_i-|\alpha|-2}\,  \big( 1+|s|^{1/2}r\big)^{-1-\mu_2-|\alpha|+\varepsilon}
\]
for $|\alpha|\le 1$. Here $\varepsilon$ is an arbitrarily small positive number, $\varepsilon=0$ if $d_i=0$.
If $\lambda_i=2m+1$, then the same estimates with $\mu_1=0$ instead of $\mu_2$ are valid.
\end{Co}

The last corollary enables us to improve the estimate (\ref{4At1}) for the coefficients $c_{i,k}$ and the estimate (\ref{3At1})
for the remainder $(v,q)$ in Theorem \ref{At1} if $\gamma>-\mu_2-\frac 12$. As in Theorem \ref{At1}, let $I_{\gamma}$ denote
the set of all $i$ such that $0 < \mbox{Re}\, \lambda_i < \frac 12 - \gamma$.

\begin{Th} \label{At2}
Let $(u,p) \in E_\beta^2(K) \times V_\beta^1(K)$ be a solution of the problem {\em (\ref{par3})} with the right-hand sides {\em (\ref{fandg})}.
We assume that $\mbox{\em Re}\, s \ge 0$, $s\not= 0$, and that there are no eigenvalues of the pencil ${\cal L}(\lambda)$ on the line
$\mbox{\em Re}\, \lambda = \frac 12 -\gamma$.

{\em 1)} If $\lambda_1<1$ and $-\frac 12 < \gamma <\frac 12-\lambda_1 < \beta < \lambda_1 +\frac 32$, then
\begin{equation} \label{1At2}
(u,p) = \zeta_s \sum_{i\in I_{\gamma}} \sum_{k=1}^{\kappa_i} c_{i,k}(s)\, \big( u_{i,k}^{(0)},p_{i,k}^{(0)}\big) + (v,q),
\end{equation}
where $v\in E_\gamma^2(K)$, $q\in V_\gamma^1(K)$, and
\begin{equation} \label{4At2}
\| v\|^2_{V_\gamma^2(K)} + |s|^2\, \| v\|^2_{V_\gamma^0(K)} + \| q\|^2_{V_\gamma^1(K)}
   \le c\, \Big( \| f\|^2_{V_\gamma^0(K)} +  \| g\|^2_{V_\gamma^1(K)}+ |s|^2\, \| g\|^2_{(V_{-\gamma}^1(K))^*}\Big)
\end{equation}
with a constant $c$ independent of $f$, $g$ and $s$.

{\em 2)} If $-\mu_2-\frac 12 < \gamma <-\frac 12 \le \frac 12-\lambda_1 < \beta < \lambda_1+\frac 32$, the eigenvalue $\lambda=1$ is simple and there
are no odd integer eigenvalues of the pencil ${\cal L}(\lambda)$ in the interval $1 < \lambda <\frac 12-\gamma$, then
\begin{equation} \label{1At2a}
u = \zeta_s \sum_{\substack{i\in I_{\gamma} \\ \lambda_i\not=1}} \ \sum_{k=1}^{\kappa_i} c_{i,k}(s)\,  u_{i,k,\gamma} + v,
\end{equation}
where  $v\in E_\gamma^2(K)$ and
\begin{equation} \label{2At2}
\| v\|^2_{V_\gamma^2(K)} + |s|^2\, \| v\|^2_{V_\gamma^0(K)}
   \le c\, \Big( \| f\|^2_{V_\gamma^0(K)} +  \| g\|^2_{V_\gamma^1(K)}+ |s|^2\, \| g\|^2_{(V_{-\gamma}^1(K))^*}\Big)
\end{equation}
with a constant $c$ independent of $f$, $g$ and $s$.

The coefficients $c_{i,k}$ in {\em (\ref{1At2})} and {\em (\ref{1At2a})} satisfy the estimate
\begin{equation} \label{3At2}
|c_{i,k}(s)|^2 \le  c\, |s|^{\gamma+ \mbox{\em \scriptsize Re}\,\lambda_i -1/2}\, \big( 1+|\log|s||\big)^{2d_i}\, \Big( \| f\|^2_{V_\gamma^0(K)}
  +  \| g\|^2_{V_\gamma^1(K)} + |s|^2\, \| g\|^2_{(V_{-\gamma}^1(K))^*} \Big)
\end{equation}
Here, $d_i$ is the same number as in Lemma {\em \ref{Al5}}. In particular, $d_i=0$ for $0<\lambda_i<1$.
\end{Th}

P r o o f.
By Theorem \ref{At1}, the solution $(u,p$ has the representation (\ref{1At1}).
Obviously, $N_{i,\gamma}=0$ for $i\ge 1$ and $\gamma>-\frac 32$. In this case, the pair $(u_{i,k,\gamma},p_{i,k,\gamma})$ coincides
with  $(u_{i,k}^{(0)},p_{i,k}^{(0)})$. The coefficients in (\ref{1At1}) are given by the formula
\[
c_{i,k}(s) = -\int_K \Big( f\cdot (\zeta_s v_{i,k} +v'_{i,k}) + g\, ( \zeta_s q_{i,k}+q'_{i,k})\Big)\, dx.
\]
Using the estimates
\[
\big| \zeta_s v_{i,k}(x,s)\big| \le c\, r^{-\mbox{\scriptsize Re}\, \lambda_i-1} (1+\log r)^{d_i}, \quad
\big| \zeta_s q_{i,k}(x,s)\big| \le c\, r^{-\mbox{\scriptsize Re}\, \lambda_i-2} (1+\log r)^{d_i},
\]
we easily obtain
\[
\| \zeta_s v_{i,k}\|^2_{V_{-\gamma}^0(K)} + \| \zeta_s q_{i,k}\|^2_{V_{1-\gamma}^0(K)}
  \le c\, |s|^{\gamma+\mbox{\scriptsize Re}\, \lambda_i-1/2}\, \big( 1+|\log |s||\big)^{2d_i}.
\]
for $\mbox{Re}\, \lambda_i < \frac 12 -\gamma$. Using Lemma \ref{Al5}, we obtain
\begin{eqnarray*}
\int_K r^{-2\gamma}\, |v'_{i,k}|^2\, dx & \le & c\, |s|^{\mbox{\scriptsize Re}\,\lambda_i+1} \, \big( 1+|\log |s||\big)^{2d_i} \\
&& \quad \times  \Big(
\int\limits_{|s|\, r^2<1} r^{-2\gamma}\, (|s|\, r^2)^{-\delta+1/2}\, dx + \int\limits_{|s|\, r^2>1} r^{-2\gamma}\, (|s|\, r^2)^{-2-\mu_2}\, dx\Big)\\
& \le & c\,  |s|^{\gamma + \mbox{\scriptsize Re}\,\lambda_i-1/2} \, \big( 1+|\log |s||\big)^{2d_i}
\end{eqnarray*}
if $\lambda_i\not=2m+1$ for integer $m\ge 0$ and $\gamma>-\mu_2-\frac 12$ (since $-2\gamma-2\delta+1>-2\gamma>-3$). Furthermore,
\begin{eqnarray*}
\| \zeta_s q'_{i,k}\|^2_{V_{1-\gamma}^0(K)} & \le & 
 c\, |s|^{\mbox{\scriptsize Re}\,\lambda_i+2} \, \big( 1+|\log |s||\big)^{2d_i}\int\limits_{|s|r^2<1} r^{2-2\gamma}\, (|s| \, r^2)^{-\delta-1/2}\, dx \\
& \le & c\, |s|^{\gamma+\mbox{\scriptsize Re}\,\lambda_i-1/2} \, \big( 1+|\log |s||\big)^{2d_i}
\end{eqnarray*}
and
\begin{eqnarray*}
\| (1-\zeta_s) q'_{i,k}\|^2_{V_{-\gamma}^1(K)} & \le &
 c\, \int\limits_{|s|r^2>1/2} \big( r^{-2-2\gamma}\, |q'_{i,k}|^2 + r^{-2\gamma} |\nabla q'_{i,k}|^2\big)\, dx \\
& \le & c\, |s|^{\gamma+\mbox{\scriptsize Re}\,\lambda_i+3/2} \, \big( 1+|\log |s||\big)^{2d_i}.
\end{eqnarray*}
Consequently,
\begin{eqnarray*}
&& \hspace{-2em}\Big| \int_K \big( f\cdot (\zeta_s v_{i,k}+v'_{i,k}) + g\, (\zeta_s q_{i,k}+q'_{i,k})\big)\, dx\Big|  \\
&& \le \| f\|_{V_\gamma^0(K)}\, \| \zeta_s v_{i,k}+v'_{i,k}\|_{V_{-\gamma}^0(K)} + \| g\|_{V_{\gamma-1}^0(K)}\, \| \zeta_s (q_{i,k}+q'_{i,k})\|_{V_{1-\gamma}^0(K)} \\
&& \hspace{2em}  + \ \| g\|_{(V_{-\gamma}^1(K))^*}\, \| (1-\zeta_s)q'_{i,k}\|_{V_{-\gamma}^1(K)}\\
&& \le c\, |s|^{(\gamma+\mbox{\scriptsize Re}\,\lambda_i-1/2)/2} \, \big( 1+|\log |s||\big)^{d_i}\, \Big(
  \| f\|_{V_\gamma^0(K)} + \| g\|_{V_\gamma^1(K)} + |s|\, \| g\|_{(V_{-\gamma}^1(K))^*}\Big).
\end{eqnarray*}
This proves (\ref{3At2}).

We prove (\ref{4At2}) for $\gamma>-\frac 12$. By \cite[Theorem 5.3.2]{kmr-01}, there do not exist generalized eigenvectors
to the eigenvalues $\lambda_i$, $i\in I_\gamma$. This means that $u_{i,k,\gamma}=u_{i,k}^{(0)}$ and $p_{i,k,\gamma}=p_{i,k}^{(0)}$ do not contain
logarithmic terms. The remainder $(v,q)$ is a solution of the problem
\[
(s-\Delta)\, v+ \nabla q =f', \ \ -\nabla\cdot v = g' \ \mbox{ in }K, \ \ v=0\ \mbox{ on }\partial K\backslash \{ 0\},
\]
where
\[
f' = f -  \sum_{i\in I_\gamma} \sum_{k=1}^{\kappa_i} c_{i,k}(s)\, \big( (s-\Delta)\, (\zeta_s u_{i,k}^{(0)}) + \nabla\big(\zeta_s p_{i,k}^{(0)}\big)
\]
and
\[
g' = g+  \sum_{i\in I_\gamma} \sum_{k=1}^{\kappa_i} c_{i,k}(s)\, \nabla\cdot (\zeta_s u_{i,k}^{(0)}).
\]
Using the equalities $-\Delta\, u_{i,k}^{(0)} + \nabla p_{i,k}^{(0)} = 0$ and $\nabla\cdot u_{i,k}^{(0)}=0$, we get
\begin{eqnarray*}
&& \| f'-f\|^2_{V_\gamma^0(K)} + \| g'-g\|^2_{V_\gamma^1(K)} + |s|^2\, \| G-g\|^2_{(V_{-\gamma}^1(K))^*} \\
&&  \le c\, \sum_{i\in J_\gamma} \sum_{k=1}^{\kappa_i} \big| c_{i,k}(s,\tau)\big|^2\, |s|^{-\gamma-\mbox{\scriptsize Re}\,\lambda_i+1/2}.
\end{eqnarray*}
This together with (\ref{3At2}) yields
\[
\| f'\|^2_{V_\gamma^0(K)} + \| g'\|^2_{V_\gamma^1(K)} + |s|^2\, \| g'\|^2_{(V_{-\gamma}^1(K))^*}
 \le c\, \Big(  \| f\|^2_{V_\gamma^0(K)} + \| g\|^2_{V_\gamma^1(K)} +   |s|^2\, \| g\|^2_{(V_{-\gamma}^1(K))^*} \Big).
\]
Since the operator $A_\gamma$ is injective (see Theorem \ref{t1}) and has closed range, we obtain the estimate (\ref{4At2}) for the remainder $(v,q)$.

We consider the case $\gamma<-1/2$. Then the strip $0 < \mbox{Re}\, \lambda < \frac 12 -\gamma$ contains the
simple eigenvalue $\lambda_j=1$ of the pencil ${\cal L}(\lambda)$. Since $u_{j,1,\gamma}(x,s)=u_{j,1}^{(0)}(x)=0$ and $p_{j,1,\gamma}(x,s)=p_{j,1}^{(0)}(x)=1$
(see Remark \ref{Ar1}), it follows from Theorem \ref{At1} that
\[
u(x) = \zeta_s(x) \sum_{i\in I_\gamma\backslash\{j\}} \sum_{k=1}^{\kappa_i} c_{i,k}(s)\, u_{i,k,\gamma}(x,s)+ v(x)
\]
and
\[
p(x) = \zeta_s(x) \sum_{i\in I_\gamma\backslash\{j\}} \sum_{k=1}^{\kappa_i} c_{i,k}(s)\, p_{i,k,\gamma}(x,s)+ c_{j,1}(s) +q'(x),
\]
where $v\in E_\gamma^2(K)$ and $q'=q+(\zeta_s-1)\, c_{j,1}(s) \in V_\gamma^1(K)$. By Remark \ref{r1}, we can write these decompositions of $u$ and $p$ also in
the form
\begin{eqnarray*}
u(x) & = & \zeta_s(x) \sum_{i\in I_\gamma\backslash\{j\}} \sum_{k=1}^{\kappa_i} C_{i,k}(s,\tau)\, U_{i,k,\gamma}(x,s,\tau)+ v(x), \\
p(x) & = & \zeta_s(x) \sum_{i\in I_\gamma\backslash\{j\}} \sum_{k=1}^{\kappa_i} C_{i,k}(s,\tau)\, P_{i,k,\gamma}(x,s)+ c_{j,1}(s) +q'(x),
\end{eqnarray*}
where $\tau$ is an arbitrary positive real number, $v$ and $q'$ are independent of $\tau$. The coefficients $C_{i,k}(s,\tau)$
are given by the formula
\[
C_{i,k}(s,\tau) = -\int_K \big( f(x)\cdot V^*_{i,k}(x,s,\tau) + g(x)\, Q^*_{i,k}(x,s,\tau)\big)\, dx,
\]
where $V^*_{i,k}(x,s,\tau)= \tau^{\lambda_i+1}\, v^*_{i,k}(\tau x,\tau^{-2}s)$ and $Q^*_{i,k}(x,s,\tau)= \tau^{\lambda_i+2}\, q^*_{i,k}(\tau x,\tau^{-2}s)$
(see Remark \ref{Ar3}). One easily checks that
\[
\| V_{i,k}^*(\cdot,s,\tau)\|^2_{V_{-\gamma}^0(K)} = \tau^{2(\gamma+\mbox{\scriptsize Re}\, \lambda_i)-1}\int_K r^{-2\gamma}\, \big|
  v_{i,k}^*(x,\tau^{-2}s)\big|^2\, dx.
\]
Using the estimates for $v'_{i,k}$ and $\zeta_s v_{i,k}$ given above, we obtain
\[
\| V_{i,k}^*(\cdot,s,\tau)\|^2_{V_{-\gamma}^0(K)} \le c\, |s|^{\gamma + \mbox{\scriptsize Re}\,\lambda_i-1/2} \, \Big( 1+\Big|\log \frac{|s|}{\tau^2}\Big|\Big)^{2d_i}.
\]
Analogously,
\[
\| \zeta_s Q^*_{i,k}\|^2_{V_{1-\gamma}^0(K)} + |s|^{-2}\, \| (1-\zeta_s) Q^*_{i,k}\|^2_{V_{-\gamma}^1(K)}
  \le c\, |s|^{\gamma + \mbox{\scriptsize Re}\,\lambda_i-1/2} \, \Big( 1+\Big|\log \frac{|s|}{\tau^2}\Big|\Big)^{2d_i}.
\]
Hence, in the case $\tau=\sqrt{|s|}$, we obtain
\[
\big| C_{i,k}(s,\tau)\big|^2 \le  c\, |s|^{\gamma + \mbox{\scriptsize Re}\,\lambda_i-1/2} \, 
  \Big(  \| f\|^2_{V_\gamma^0(K)} + \| g\|^2_{V_\gamma^1(K)} + |s|^2\, \| g\|^2_{(V_{-\gamma}^1(K))^*}\Big)
\]
with a constant $c$ independent of $f$, $g$ and $s$.

The pair $(v,q')$ is a solution of the problem
\[
(s-\Delta)\, v+ \nabla q' =F, \ \ -\nabla\cdot v = G \ \mbox{ in }K, \ \ v=0\ \mbox{ on }\partial K\backslash \{ 0\},
\]
where
\[
F = f -  \sum_{i\in I_\gamma\backslash\{j\}} \sum_{k=1}^{\kappa_i} C_{i,k}(s,\tau)\, \big( (s-\Delta)\, (\zeta_s U_{i,k,\gamma} + \nabla(\zeta_s P_{i,k,\gamma})\big)
\]
and
\[
G = g+  \sum_{i\in I_\gamma\backslash\{j\}} \sum_{k=1}^{\kappa_i} C_{i,k}(s,\tau)\, \nabla\cdot (\zeta_s U_{i,k,\gamma}).
\]
Using the equalities
\[
(s-\Delta)\, U_{i,k,\gamma} + \nabla P_{i,k,\gamma} = s\, (sr^2)^{N_{i,\gamma}} U_{i,k}^{N_{i,\gamma}}, \quad \nabla\cdot U_{i,k,\gamma}
\]
(see Lemma \ref{Al1}), we obtain the estimate
\begin{eqnarray*}
&& \| F-f\|^2_{V_\gamma^0(K)} + \| G-g\|^2_{V_\gamma^1(K)} + |s|^2\, \| G-g\|^2_{(V_{-\gamma}^1(K))^*} \\
&&  \le c\, \sum_{i\in J_\gamma\backslash\{j\}} \sum_{k=1}^{\kappa_i} \big| C_{i,k}(s,\tau)\big|^2\, |s|^{-\gamma-\mbox{\scriptsize Re}\,\lambda_i+1/2}\,
  \Big( 1+\Big|\log \frac{|s|}{\tau^2}\Big|\Big)^{2d_i}
\end{eqnarray*}
since $\gamma+\mbox{Re}\, \lambda_i +2N_{i,\gamma}>-\frac 32$. Thus, in the case $\tau=\sqrt{|s|}$, we obtain
\[
\| F\|^2_{V_\gamma^0(K)} + \| G\|^2_{V_\gamma^1(K)} + |s|^2\, \| G\|^2_{(V_{-\gamma}^1(K))^*}
 \le c\, \Big(  \| f\|^2_{V_\gamma^0(K)} + \| g\|^2_{V_\gamma^1(K)} +   |s|^2\, \| g\|^2_{(V_{-\gamma}^1(K))^*} \Big).
\]
Under the conditions of the theorem, the operator $A_\gamma$ is injective (see Theorem \ref{t1}) and has closed range
(see \cite[Theorem 2.1]{k/r-16}). This implies
\[ 
\| v\|^2_{V_\gamma^2(K)} + |s|^2\, \| v\|^2_{V_\gamma^0(K)} + \| q'\|^2_{V_\gamma^1(K)}
   \le c\, \Big( \| f\|^2_{V_\gamma^0(K)} +  \| g\|^2_{V_\gamma^1(K)}+ |s|^2\, \| g\|^2_{(V_{-\gamma}^1(K))^*}\Big).
\]
In particular, the estimate (\ref{2At2}) holds. The proof of the theorem is complete. \hfill $\Box$

\begin{Rem} \label{Ar4}
{\em 1) Conditions for the simplicity of the eigenvalue $\lambda=1$ can be found in \cite[Section 5.5]{kmr-01}. For example, this eigenvalue
is simple if $\bar{\Omega}$ is a subset of a halfsphere.

2) By \cite[Theorem 5.3.2]{kmr-01}, the eigenvalues of the pencil ${\cal L}(\lambda)$ in the strip $0<\mbox{Re}\, \lambda <1$ do not have
generalized eigenvectors. Thus, the singular functions in (\ref{1At2}) have the form
\[
u_{i,k}^{(0)}(x)= r^{\lambda_i}\, \phi_{i,k}(\omega), \quad p_{i,k}^{(0)}(x)= r^{\lambda_i-1}\, \psi_{i,k}(\omega)
\]
where $(\phi_{i,k},\psi_{i,k})$ are eigenvectors of the pencil ${\cal L}(\lambda)$ corresponding to the eigenvalue
$\lambda_i$.}
\end{Rem}

\section{The time-dependent problem}

Finally, we consider the problem problem (\ref{stokes1}), (\ref{stokes2}).
Using the results of the last section, we obtain the asymptotics of the solution near the vertex of the cone.

\subsection{Asymptotics of the solution near the vertex of the cone}

Let $Q=K\times {\Bbb R}_+ = K\times (0,\infty)$.
We denote by $W_\beta^{2,1}(Q)$ the weighted Sobolev space of all functions $u=u(x,t)$ on $Q$ with finite norm
\[
\| u\|_{W_\beta^{2,1}(Q)} = \Big( \int_0^\infty  \big( \| u(\cdot,t) \|^2_{V_\beta^2(K)}+\| \partial_t u(\cdot,t) \|^2_{V_\beta^0(K)}\big)\, dt \Big)^{1/2}.
\]
Note that $u(\cdot,0) \in V_\beta^1(K)$ for $u\in W_\beta^{2,1}(Q)$ (see \cite[Proposition 3.1]{kozlov-89}).
The space $\stackrel{\circ}{W}\!{}_\beta^{2,1}(Q)$ is the subspace of all $u\in W_\beta^{2,1}(Q)$ satisfying
the condition $u(x,0)=0$ for $x\in K$. By \cite[Proposition 3.4]{kozlov-89}, the Laplace transform realizes an isomorphism from
$\stackrel{\circ}{W}\!{}_\beta^{2,1}(Q)$ onto the space $H_\beta^{2}$ of all holomorphic functions $\tilde{u}(x,s)$ for $\mbox{Re}\, s >0$ with values
in $E_\beta^2(K)$ and finite norm
\[
\| \tilde{u}\|_{H_\beta^{2}} = \sup_{\gamma>0} \Big( \int_{-\infty}^{+\infty} \big(
  \| \tilde{u}(\cdot,\gamma+i\tau)\|^2_{V_\beta^2(K)} +|\gamma+i\tau|^2\, \| \tilde{u}(\cdot,\gamma+i\tau)\|^2_{V_\beta^0(K)}\big) \, d\tau  \Big)^{1/2}.
\]
The proof of the analogous result in nonweighted spaces can be found in \cite[Theorem 8.1]{av-64}.

The following existence and uniqueness theorem  was proved in \cite[Theorem 3.1]{k/r-16}.

\begin{Th} \label{t6}
Suppose that $f\in L_2\big( {\Bbb R}_+,V_\beta^0(K)\big)$, $g\in L_2\big({\Bbb R}_+,V_\beta^1(K)\big)$,
$\partial_t g\in L_2\big({\Bbb R}_+,(V_{-\beta}^1(K))^*\big)$ and $g(x,0)=0$ for $x\in K$, where
$\beta$ satisfies the inequalities {\em (\ref{interval})}.
In the case $\beta>1/2$, we assume in addition that $\int_K g(x,t)\, dx =0$ for almost all $t$.
Then there exists a uniquely determined solution $(u,p) \in W_\beta^{2,1}(Q) \times L_2\big({\Bbb R}_+, V_\beta^1(K)\big)$
of the problem {\em (\ref{stokes1}), (\ref{stokes2})} satisfying the estimate
\[ 
\| u\|_{W_\beta^{2,1}(Q)} + \| p\|_{L_2({\Bbb R}_+,V_\beta^1(K))}  \le c\, \Big( \| f\|_{W_\beta^{0,0}(Q)}
  + \| g \|_{L_2({\Bbb R}_+,V_\beta^1(K))} +  \| \partial_t g\|_{L_2({\Bbb R}_+,(V_{-\beta}^1(K))^*)}\Big)
\]
with a constant $c$ independent of $f$ and $g$.
\end{Th}

The condition on $\beta$ in Theorem \ref{t6} can be weakened if $\lambda_1=1$ and $\lambda_1$ is a simple eigenvalue of the pencil ${\cal L}(\lambda)$
(see \cite{k/r-18a}).

We consider the solution $(u,p)\in  W_\beta^{2,1}(Q) \times L_2\big({\Bbb R}_+, V_\beta^1(K)\big)$ of the problem (\ref{stokes1}), (\ref{stokes2}), where
\begin{equation} \label{E1}
f\in L_2\big( {\Bbb R}_+,V_\beta^0(K) \cap V_\gamma^0(K)\big), \quad
g\in L_2\big({\Bbb R}_+,V_\beta^1(K) \cap V_\gamma^1(K)\big)
\end{equation}
and
\begin{equation} \label{E2}
\partial_t g\in L_2\big({\Bbb R}_+,(V_{-\beta}^1(K))^* \cap (V_{-\gamma}^1(K))^*\big),
\end{equation}
$\beta$ satisfies the inequalities (\ref{interval}) and $\gamma < \beta$. If $\gamma >\frac 12 - \lambda_1$, then it follows from \cite[Theorem 3.2]{k/r-16} that
$u\in W_\gamma^{2,1}(Q)$ and $p\in L_2\big({\Bbb R}_+, V_\gamma^1(K)\big)$. For $\gamma <\frac 12 - \lambda_1$, we can use the decomposition in
Theorem \ref{At2} for the Laplace transforms of $u$ and $p$.

Let $\psi$ be a $C^\infty$-function on $(-\infty,+\infty)$ with support in the interval $[0,1]$ satisfying the conditions
\[
\int_0^1 \psi(t)\, dt = 1, \quad \int_0^1 t^j\, \psi(t)\, dt = 0 \ \mbox{ for }j=1,\ldots,N,
\]
where $N$ is an integer, $2N> \frac 12 - \gamma-\lambda_1$. By $\tilde{\psi}$, we denote the Laplace transform of $\psi$.
The function $\tilde{\psi}$ is analytic in ${\Bbb C}$ and satisfies the conditions
$\tilde{\psi}(0)=1$, $\tilde{\psi}^{(j)}(0)=0$ for $j=1,\ldots,N$. Since $s^n \tilde{\psi}^{(j)}(s)$ is the Laplace transform of the function
$(-1)^j\, \frac{d^n}{dt^n}\big( t^j\, \psi(t)\big)$, it follows that
\[
\big|\tilde{\psi}^{(j)}(s)\big| \le c_{j,n}\, |s|^{-n}
\]
for every $j\ge 0$ and $n\ge 0$, where $c_{j,n}$ is independent of $s$, $\mbox{Re}\, s\ge 0$. We consider the functions
\[
\tilde{K}_{i,k}(x,y,s) =\tilde{\psi}(sr^2)\, v^*_{i,k}(y,s) \quad \mbox{and}\quad \tilde{T}_{i,k}(x,y,s) =\tilde{\psi}(sr^2)\, q^*_{i,k}(y,s),
\]
where $(v^*_{i,k},q^*_{i,k})$ is the special solution of the problem (\ref{hom}) introduced in Section 2.2.
The functions $\tilde{K}_{i,k}(x,y,s)$ and $\tilde{T}_{i,k}(x,y,s)$  are the Laplace transforms of
\begin{equation} \label{2l2}
K_{i,k}(x,y,t) = \frac{1}{2\pi i} \int_{-i\infty}^{+i\infty} e^{st}\, \tilde{K}_{i,k}(x,y,s)\, ds \ \mbox{ and } \
T_{i,k}(x,y,t) = \frac{1}{2\pi i} \int_{-i\infty}^{+i\infty} e^{st}\, \tilde{K}_{i,k}(x,y,s)\, ds,
\end{equation}
respectively. Obviously, $K_{i,k}(x,y,t)$ and $H_{i,k}(x,y,t)$ depend only on $r=|x|$, $y$ and $t$.

Finally, let $u_{i,k,\gamma}(x,s)$, $p_{i,k,\gamma}(x,s)$ be the functions defined in (\ref{Uikg}). These functions are polynomials in $s$
of degree $N_{i,\gamma}$. If we replace $s$ by $\partial_t$ we get the differential operators
$u_{i,k,\gamma}(x,\partial_t)$ and $p_{i,k,\gamma}(x,\partial_t)$, In the case $\gamma+\mbox{Re}\, \lambda_i>-\frac 32$, we have
$u_{i,k,\gamma}(x,\partial_t)=u_{i,k}^{(0)}(x)$ and $p_{i,k,\gamma}(x,\partial_t)=p_{i,k}^{(0)}(x)$.

As a consequence of Theorem \ref{At2}, we obtain the following assertion.

\begin{Th} \label{t8}
Let $(u,p)\in  W_\beta^{2,1}(Q) \times L_2\big({\Bbb R}_+, V_\beta^1(K)\big)$ be a solution of the problem {\em (\ref{stokes1}), (\ref{stokes2})}
with data $f,g$ satisfying {\em (\ref{E1}), (\ref{E2})}. We assume that $\frac 12-\lambda_1 < \beta < \lambda_1+\frac 32$, $-\mu_2-\frac 12 <\gamma<\beta$ and that
the line $\mbox{\em Re}\, \lambda = \frac 12 -\gamma$ is free of eigenvalues of the pencil ${\cal L}(\lambda)$. By $I_\gamma$, we denote the set of all
$i$ such that $0<\mbox{\em Re}\, \lambda_i < \frac 12 -\gamma$.

{\em 1)} If $\lambda_1<1$ and $-\frac 12 < \gamma < \beta$, then
\begin{equation} \label{8t8}
(u,p) =  \sum_{i\in I_\gamma} \sum_{k=1}^{\kappa_i} H_{i,k}(r,t) \,  \big( u_{i,k}^{(0)}(x),p_{i,k}^{(0)}(x)\big)   + (v,q),
\end{equation}
where $H_{i,k}$ is given by {\em (\ref{7t8})}, $v\in  W_\gamma^{2,1}(Q)$ and $q\in L_2\big({\Bbb R}_+, V_\gamma^1(K)\big)$.

{\em 2)} If $-\mu_2-\frac 12 < \gamma <-\frac 12$, the eigenvalue $\lambda=1$ is simple and there
are no odd integer eigenvalues of the pencil ${\cal L}(\lambda)$ in the interval $1 < \lambda <\frac 12-\gamma$, then
\begin{equation} \label{1t8}
u(x,t) = \sum_{\substack{i\in I_\gamma \\ \lambda_i\not=1}} \ \sum_{k=1}^{\kappa_i} u_{i,k,\gamma}(x,\partial_t)\, H_{i,k}(r,t) + v(x,t),
\end{equation}
where $H_{i,k}$ is given by {\em (\ref{7t8})} and $v\in W_\gamma^{2,1}(Q)$.

The functions $v$ and $q$ in {\em (\ref{8t8})} and {\em (\ref{1t8})} satisfy the estimates
\begin{eqnarray} \label{4t8}
\| v\|_{W_\gamma^{2,1}(Q)} \le c\, \Big( \| f\|_{L_2({\Bbb R}_+,V_\gamma^0(K))}
  + \| g\|_{L_2({\Bbb R}_+, V_\gamma^1(K))} + \| \partial_t g\|_{L_2\big({\Bbb R}_+, (V_{-\gamma}^1(K))^*\big)}\Big)
\end{eqnarray}
and
\[
\| q\|_{L_2({\Bbb R}_+, V_\gamma^1(K))} \le c\, \Big( \| f\|_{L_2({\Bbb R}_+,V_\gamma^0(K))}
  + \| g\|_{L_2({\Bbb R}_+, V_\gamma^1(K))} + \| \partial_t g\|_{L_2\big({\Bbb R}_+, (V_{-\gamma}^1(K))^*\big)}\Big)
\]
with a constant $c$ independent of $f$ and $g$.
\end{Th}

P r o o f.
We may restrict ourselves to the case $\gamma<\frac 12-\lambda_1$. If $\gamma>\frac 12-\lambda_1$, then the set $I_\gamma$ is empty,
and the assertion of the theorem follows immediately from \cite[Theorem 3.2]{k/r-16}.
By Theorem \ref{At2}, the Laplace transform of $u$ admits the decomposition
\begin{equation} \label{3t8}
\tilde{u}(x,s) = \zeta(|s|r^2) \sum_{i\in I_{\gamma}} \sum_{k=1}^{\kappa_i} c_{i,k}(s)\,
   u_{i,k,\gamma}(x,s) + V(x,s),
\end{equation}
where $V$ satisfies the estimate (\ref{2At2}) with $\tilde{f}$, $\tilde{g}$ instead of $f$ and $g$ and
\begin{equation} \label{2t8}
c_{i,k}(s) = -\int_K \big( \tilde{f}(y,s)\cdot v^*_{i,k}(y,s) + \tilde{g}(y,s)\, q^*_{i,k}(y,s)\big)\, dy.
\end{equation}
Here $u_{i,k,\gamma}=0$ for $\lambda_i=1$.
As in the proof of Theorem \ref{At2}, we can write (\ref{3t8}) in the form
\begin{equation} \label{6t8}
\tilde{u}(x,s) = \zeta(|s|r^2) \sum_{i\in I_{\gamma}} \sum_{k=1}^{\kappa_i} C_{i,k}(s,\tau)\,
   U_{i,k,\gamma}(x,s,\tau) + V(x,s),
\end{equation}
with arbitrary positive $\tau$ and with the same $V$ as in (\ref{3t8}). For $\tau=|s|^{1/2}$, the function
$C_{i,k}(s,\tau)$ satisfies the estimate
\[
\big| C_{i,k}(s,\tau)\big|^2 \le  c\, |s|^{\gamma + \mbox{\scriptsize Re}\,\lambda_i-1/2} \,
  \Big(  \| \tilde{f}(\cdot,s)\|^2_{V_\gamma^0(K)} + \| \tilde{g}(\cdot,s)\|^2_{V_\gamma^1(K)} + |s|^2\, \| \tilde{g}(\cdot,s)\|^2_{(V_{-\gamma}^1(K))^*}\Big).
\]
Using the inequalities
\[
\big| \partial_r^j \big( \zeta(|s|r^2)-\tilde{\psi}(sr^2)\big)\big| \le c\, r^{-j} \, |sr^2|^{N+1} \ \mbox{for }|sr^2|<1
\]
and
\[
\big| \partial_r^j \big( \zeta(|s|r^2)-\tilde{\psi}(sr^2)\big)\big| \le c\, r^{-j} \, |sr^2|^{-n} \ \mbox{for }|sr^2|>1,
\]
we obtain the estimate
\begin{eqnarray*}
&& \big\| \big( \zeta(|s|r^2)- \tilde{\psi}(sr^2)\big)\, U_{i,k,\gamma}(\cdot,s,\tau)\big\|^2_{V_\gamma^2(K)}
+ |s|^2\, \big\| \big( \zeta(|s|r^2)- \tilde{\psi}(sr^2)\big)\,U_{i,k,\gamma}(\cdot,s,\tau)\big\|^2_{V_\gamma^0(K)} \\
&& \le c\, |s|^{-\gamma-\mbox{\scriptsize Re}\, \lambda_i+1/2}
\end{eqnarray*}
for $\tau=|s|^{1/2}$. Hence, it follows from (\ref{6t8}) that
\[
\tilde{u}(x,s) = \tilde{\psi}(sr^2) \sum_{i\in I_{\gamma}} \sum_{k=1}^{\kappa_i} C_{i,k}(s,|s|^{1/2})\,
   U_{i,k,\gamma}(x,s,|s|^{1/2}) + \tilde{v}(x,s),
\]
or, equivalently,
\[
\tilde{u}(x,s) = \tilde{\psi}(sr^2) \sum_{i\in I_{\gamma}} \sum_{k=1}^{\kappa_i} c_{i,k}(s)\, u_{i,k,\gamma}(x,s) + \tilde{v}(x,s),
\]
where
\[
\| \tilde{v}(\cdot,s)\|^2_{V_\gamma^2(K)} + |s|^2\, \| \tilde{v}(\cdot,s)\|^2_{V_\gamma^0(K)}
   \le c\, \Big( \| \tilde{f}(\cdot,s)\|^2_{V_\gamma^0(K)} +  \| \tilde{q}(\cdot,s)\|^2_{V_\gamma^1(K)}
   + |s|^2\, \| \tilde{g}(\cdot,s)\|^2_{(V_{-\gamma}^1(K))^*}\Big)
\]
with a constant $c$ independent of $s$ for $\mbox{Re}\, s\ge 0$. The function $\tilde{\psi}(sr^2) c_{i,k}(s)\, u_{i,k,\gamma}(x,s)$
is the Laplace transform of $u_{i,k,\gamma}(x,\partial_t)\, H_{i,k}(r,t)$, while $\tilde{v}(x,s)$ is the Laplace transform
of a function $v\in W_\gamma^{2,1}(Q)$ satisfying the estimate (\ref{4t8}). This proves the assertion of the theorem concerning
the function $u$. Analogously, the assertion concerning $p$ holds.  \hfill $\Box$ \\

The condition on $\gamma$ in Theorem \ref{t8} can be weakened if we apply Theorem \ref{At1} instead of Theorem \ref{At2}.
But then we need additional conditions on the $t$-derivatives of the data $f$ and $g$.

\begin{Th} \label{t8a}
Let $(u,p)\in  W_\beta^{2,1}(Q) \times L_2\big({\Bbb R}_+, V_\beta^1(K)\big)$ be the solution of the problem {\em (\ref{stokes1}), (\ref{stokes2})}
with data $f,g$ satisfying {\em (\ref{E1}), (\ref{E2})}. We assume that $\gamma < \frac 12-\lambda_1 < \beta < \frac 12$
and that $\lambda=\frac 12 - \gamma$ is not an eigenvalue of the pencil ${\cal L}(\lambda)$. Furthermore, we assume that
\[
\partial_t^n f \in L_2({\Bbb R}_+,V_{\gamma+2n}^0(K)), \ \ \partial_t^n g \in L_2({\Bbb R}_+,V_{\gamma+2n}^1(K)), \ \
\partial_t^{n+1} g \in L_2\big({\Bbb R}_+,(V_{-\gamma-2n}^0(K))^*\big),
\]
where $n$ is the smallest integer such that $\beta-\gamma\le 2n$. Then
\[ 
( u,p) = \sum_{i\in I_{\gamma}} \sum_{k=1}^{\kappa_i} \big( u_{i,k,\gamma}(x,\partial_t), p_{i,k,\gamma}(x,\partial_t)\big)\, H_{i,k}(r,t) + (v,q),
\]
where $I_{\gamma}$ is the set of all $i$ such that $0 < \lambda_i < \frac 12 - \gamma$,  $H_{i,k}$ is the function {\em (\ref{7t8})},
$v\in  W_\gamma^{2,1}(Q)$ and $q\in L_2\big({\Bbb R}_+, V_\gamma^1(K)\big)$.
\end{Th}

P r o o f.
Let $(\tilde{u},\tilde{p})$ be the Laplace transform of $(u,p)$. Then Theorem \ref{At1} implies
\[
(\tilde{u},\tilde{p}) = \sum_{i\in I_{\gamma}} \sum_{k=1}^{\kappa_i} \big( u_{i,k,\gamma}(x,s), p_{i,k,\gamma}(x,s)\big)\, \tilde{\psi}(sr^2)\, c_{i,k}(s)
  + (\tilde{v},\tilde{q}),
\]
where $c_{i,k}(s)$ is the function (\ref{2t8}) and $(\tilde{v},\tilde{q})$ satisfies the estimate (\ref{5At1})
(with $\tilde{f},\tilde{g}$ instead of $f,g$).  Using the inequality
\[
\| \tilde{f}(\cdot,s)\|^2_{V_\gamma^0(K)} + |s|^{\beta-\gamma}\, \| \tilde{f}(\cdot,s)\|^2_{V_\beta^0(K)}
\le 2\, \Big(  \| \tilde{f}(\cdot,s)\|^2_{V_\gamma^0(K)} +  \| s^n \tilde{f}(\cdot,s)\|^2_{V_{\gamma+2n}^0(K)}\Big)
\]
and analogous inequalities for the norms of $\tilde{g}$, we obtain
\begin{eqnarray*}
&& \hspace{-2em}\| v\|_{W_\gamma^{2,1}(Q)} +\| q\|_{L_2\big({\Bbb R}_+, V_\gamma^1(K)\big)}
\le c\, \Big( \| f\|_{L_2\big({\Bbb R}_+,V_\gamma^0(K)\big)}
  + \| g\|_{L_2\big({\Bbb R}_+, V_\gamma^1(K)\big)} + \| \partial_t g\|_{L_2\big({\Bbb R}_+, (V_{-\gamma}^1(K))^*\big)} \\
&& \qquad + \ \| \partial_t^n f \|_{L_2({\Bbb R}_+,V_{\gamma+2n}^0(K))} + \| \partial_t^n g \|_{L_2({\Bbb R}_+,V_{\gamma+2n}^1(K))}
  + \| \partial_t^{n+1} g \|_{L_2\big({\Bbb R}_+,(V_{-\gamma-2n}^0(K))^*\big)}\Big).
\end{eqnarray*}
This proves the theorem \hfill $\Box$

\subsection{On the coefficients in the asymptotics}

We consider the coefficients $H_{i,k}$ in the asymptotics (\ref{1t8}). First, we estimate the functions $K_{i,k}$ and $T_{i,k}$
in the representation (\ref{7t8}) of $H_{i,k}$.

\begin{Le} \label{l2}
Let $K_{i,k}$ and $T_{i,k}$ be defined by {\em (\ref{2l2})}. If $\lambda_i\not=2m+1$ with a nonnegative integer $m$, then
\begin{eqnarray*}
\big| \partial_x^\alpha \partial_y^\gamma \partial_t^j K_{i,k}(x,y,t) \big| & \le & c_n \, t^{-(\mbox{\em \scriptsize Re}\, \lambda_i+|\alpha|+|\gamma|+2j+3)/2} \,
  \Big(1+\frac{|x|}{\sqrt{t}}\Big)^{-n} \,  \Big( \frac{|y|}{|y|+\sqrt{t}} \Big)^{-\mbox{\em \scriptsize Re}\, \lambda_i-|\gamma|-1} \\
&& \times \ \Big( 1+\frac{|y|}{\sqrt{t}}\Big)^{-2-\mu_2+\varepsilon} \, \big( 1+\big|\log|y|\big|\big)^{d_i}
\end{eqnarray*}
for $|\gamma|\le 2$ and
\begin{eqnarray*}
\big| \partial_x^\alpha \partial_y^\gamma \partial_t^j T_{i,k}(x,y,t) \big| & \le & c_n \, t^{-(\mbox{\em \scriptsize Re}\, \lambda_i+|\alpha|+|\gamma|+2j+4)/2} \,
  \Big(1+\frac{|x|}{\sqrt{t}}\Big)^{-n} \,  \Big( \frac{|y|}{|y|+\sqrt{t}} \Big)^{-\mbox{\em \scriptsize Re}\, \lambda_i-|\gamma|-2} \\
&& \times \ \Big( 1+\frac{|y|}{\sqrt{t}}\Big)^{-1-\mu_2-|\gamma|+\varepsilon} \, \big( 1+\big|\log|y|\big|\big)^{d_i}
\end{eqnarray*}
for $|\gamma|\le 1$. Here $n$ is an arbitrary nonnegative number, the constant $c_n$ is independent of $x,y,t$. Furthermore,
$\varepsilon$ is an arbitrary positive number, $\varepsilon=0$ if $d_i=0$. If $\lambda_i=2m+1$, then the same estimates with $\mu_1=0$
instead of $\mu_2$ are valid.
\end{Le}

P r o o f.
First note that all theorems of this paper are not only valid for $\mbox{Re}\, s \ge 0$ but for all complex $s=\sigma e^{i\theta}$, where $\sigma>0$ and
$|\theta| \le \delta+\frac \pi 2$ with a sufficiently small positive number $\delta$. Therefore, one can replace the path of integration
in (\ref{2l2}) by the contour $\Gamma_{t,\delta} = \Gamma_{t,\delta}^{(1)} \cup  \Gamma_{t,\delta}^{(2)}$, where
\[
\Gamma_{t,\delta}^{(1)} = \{ s=t^{-1}\, e^{i\theta}: -\delta-\pi/2 < \theta<\delta+\pi/2 \} \quad\mbox{and}\quad
\Gamma_{t,\delta}^{(2)} = \{ s=\sigma e^{\pm i(\delta+\pi/2)}: \ \sigma > t^{-1} \}.
\]
Let $\lambda_i$ be not an odd positive integer. Using the inequality $\big|\partial_x^\alpha \tilde{\psi}(sr^2)\big|\le c\, |s|^{|\alpha|/2}\, |sr^2|^{-n}$
and Corollary \ref{Ac1}, we obtain
\begin{eqnarray*}
&& \hspace{-2em} \Big| \partial_x^\alpha \partial_y^\gamma \partial_t^j \int_{\Gamma_{t,\delta}^{(1)}} e^{st} \tilde{\psi}(sr^2)\, v^*_{i,k}(y,s)\, ds\Big|
\le c \, t^{-(\mbox{\scriptsize Re}\, \lambda_i+|\alpha|+|\gamma|+2j+3)/2} \,
  \Big(\frac{|x|}{\sqrt{t}}\Big)^{-n}  \\
&& \times \ \Big( \frac{|y|}{|y|+\sqrt{t}} \Big)^{-\mbox{\scriptsize Re}\, \lambda_i-|\gamma|-1} \
  \Big( 1+\frac{|y|}{\sqrt{t}}\Big)^{-2-\mu_2+\varepsilon} \, \big( 1+\big|\log|y|\big|\big)^{d_i}
\end{eqnarray*}
for $|\gamma|\le 2$. Furthermore,
\begin{eqnarray*}
&& \hspace{-1em}  \Big| \partial_x^\alpha \partial_y^\gamma \partial_t^j \int_{\Gamma_{t,\delta}^{(2)}} e^{st} \tilde{\psi}(sr^2)\, v^*_{i,k}(y,s)\, ds\Big|  \\
&& \le c\, \, \big( 1+\big|\log|y|\big|\big)^{d_i}  \int_{1/t}^\infty e^{-\sigma t\sin\delta}\,
  \sigma^{(\mbox{\scriptsize Re}\, \lambda_i+|\alpha|+|\gamma|+2j+1)/2}\, (\sigma r^2)^{-n}\\
&& \hspace{10em}\times \  \Big( \frac{\sigma^{1/2}|y|}{1+\sigma^{1/2}|y|}\Big)^{-\mbox{\scriptsize Re}\, \lambda_i-|\gamma|-1}\,
  \big( 1+\sigma^{1/2}|y|\big)^{-2-\mu_2+\varepsilon}\, d\sigma \\
&& \le c\, \big( 1+\big|\log|y|\big|\big)^{d_i} \Big( \frac{|y|}{|y|+\sqrt{t}}\Big)^{-\mbox{\scriptsize Re}\, \lambda_i-|\gamma|-1}\,
  \Big( 1+\frac{|y|}{\sqrt{t}}\Big)^{-2-\mu_2+\varepsilon} \\
&& \qquad \times  \int_{1/t}^\infty e^{-\sigma t\sin\delta}\, \sigma^{(\mbox{\scriptsize Re}\, \lambda_i+|\alpha|+|\gamma|+2j+1)/2}\, (\sigma r^2)^{-n}\, d\sigma
\end{eqnarray*}
which leads to the same estimate as for the integral over $\Gamma_{t,\delta}^{(1)}$.
This proves the desired estimate for $K_{i,k}$ if $\lambda_i\not=2m+1$. Analogously, $T_{i,k}$ can be estimated. In the case $\lambda_i=2m+1$,
one must replace $\mu_2$ by $\mu_1=0$ (see Corollary \ref{Ac1}). \hfill $\Box$ \\

We derive another representation of the coefficients $H_{i,k}(r,t)$ in Theorem \ref{t8}. The Laplace transform
of $H_{i,k}(r,\cdot)$ is
\begin{equation} \label{tildeH}
\tilde{H}_{i,k}(r,s) = - \int_K \big( \tilde{f}(y,s)\cdot \tilde{K}_{i,k}(x,y,s) + \tilde{g}(y,s)\, \tilde{T}_{i,k}(x,y,s)\big)\, dy
  =\tilde{\psi}(sr^2)\, \tilde{h}_{i,k}(s),
\end{equation}
where
\[
\tilde{h}_{i,k}(s) = -\int_K \big( \tilde{f}(y,s)\cdot v^*_{i,k}(y,s) + \tilde{g}(y,s)\, q^*_{i,k}(y,s)\big)\, dy
\]
is the function (\ref{2t8}). Since $\tilde{\psi}(0)=1$, it follows that $\tilde{H}_{i,k}(0,s)=\tilde{h}_{i,k}(s)$.
Obviously, $\tilde{h}_{i,k}(s)$ is the Laplace transform of
\[ 
h_{i,k}(t) = - \int_0^t \int_K \Big( f(y,\tau)\cdot V_{i,k}(y,t-\tau) + g(y,\tau)\, Q_{i,k}(y,t-\tau) \Big) \, dy\, d\tau,
\]
where $V_{i,k}(x,\cdot)$, $Q_{i,k}(x,\cdot)$ are the inverse Laplace transforms of $v^*_{i,k}(x,\cdot)$ and $q^*_{i,k}(x,\cdot)$, respectively,
\begin{eqnarray*}
V_{i,k}(x,t) & = & \frac{1}{2\pi i} (1+\partial_t) \int_{-i\infty}^{+i\infty} e^{st}\, (1+s)^{-1}\, v^*_{i,k}(x,s)\, ds, \\
Q_{i,k}(x,t) & = & \frac{1}{2\pi i} (1+\partial_t)^2 \int_{-i\infty}^{+i\infty} e^{st}\, (1+s)^{-2}\, q^*_{i,k}(x,s)\, ds
\end{eqnarray*}
(by Corollary \ref{Ac1}, the above integrals are absolutely convergent). From (\ref{tildeH}) it follows that
\[
H_{i,k} = {\cal E} h_{i,k},
\]
where the extension operator ${\cal E}$ is defined as
\[
({\cal E}h)(r,t) = r^{-2} \int_0^t \psi\Big( \frac{t-\tau}{r^2}\Big)\, h(\tau)\, d\tau.
\]

\end{document}